\documentclass[reqno]{amsart}%
\usepackage{enumerate}
\usepackage{amsfonts}
\usepackage{amsmath}
\usepackage{amssymb}
\usepackage{cases}
\usepackage{hyperref}
\usepackage{float}
\usepackage{graphicx}%
\newtheorem{theorem}{Theorem}[section]
\newtheorem{lemma}[theorem]{Lemma}
\newtheorem{proposition}[theorem]{Proposition}

\theoremstyle{definition}

\newtheorem{remark}[theorem]{Remark}

\newcommand{\norm}[1]{\left\Vert#1\right\Vert}

\numberwithin{equation}{section}
\def\fin { \vskip 0pt \hfill $\diamond$ \vskip 12pt}
\hypersetup{
colorlinks,
citecolor=blue,
filecolor=black,
linkcolor=blue,
urlcolor=magenta
}
\hypersetup{linktocpage}
\begin{document}
\title[Wave equations on hyperbolic spaces]{Well-posedness and scattering for wave equations on\\hyperbolic spaces with singular data}

\author[\ P.T. Xuan]{\ \ Pham T. Xuan}
\address{Pham Truong Xuan \hfill\break (Corresponding author) Thang Long Institute of Mathematics and Applied Sciences (TIMAS), Thang Long University, Nghiem Xuan Yem, Hanoi, Vietnam}
\author[L.C.F. Ferreira \ ]{Lucas C. F. Ferreira \ \ }
\address{Lucas C. F. Ferreira \hfill\break University of Campinas, IMECC-Department of
Mathematics, Rua S\'{e}rgio Buarque de Holanda, 651, CEP 13083-859,
Campinas-SP, Brazil} 
\email{phamtruongxuan.k5@gmail.com or xuanpt@thanglong.edu.vn}
\email{lcff@ime.unicamp.br}

\thanks{{\bf Acknowledgment.} 
L.C.F. Ferreira was supported by FAPESP (grant: 2020/05618-6) and CNPq (grant: 308799/2019-4), Brazil. }

\begin{abstract}
We consider the wave and Klein-Gordon equations on the real hyperbolic space
$\mathbb{H}^{n}$ ($n \geq2$) in a framework based on weak-$L^{p}$ spaces.
First, we establish dispersive estimates on Lorentz spaces in the context of
$\mathbb{H}^{n}$. Then, employing those estimates, we prove global well-posedness of solutions
and an exponential asymptotic stability property. Moreover, we develop a scattering theory and construct wave operators in such singular framework.\vspace{0.2cm}

 \vspace{0.2cm}


\end{abstract}
\maketitle


\font\nho=cmr10





\section{Introduction}

In the present paper, we are concerned with the nonlinear Klein-Gordon
equation on the real hyperbolic space $\mathbb{H}^{n}$ $(n\geq2)$:
\begin{equation}%
\begin{cases}
\partial_{t}^{2}u(t,x)-\Delta_{x}u(t,x)+cu(t,x)=F(u(t,x)),\\
u(0,x)=u_{0}(x),\,\,\partial_{t}u(t,x)=u_{1}(x),
\end{cases}
\label{OKlein}%
\end{equation}
where $\Delta_{x}$ stands for the Laplace-Beltrami operator associated with
the hyperbolic metric, the constant $c\geq-\dfrac{(n-1)^{2}}{4}$ and the
nonlinearity $F(u)$ satisfies
\begin{equation}
|F(0)|=0\text{ and }|F(u)-F(v)|\leq C(|u|^{b-1}+|v|^{b-1})|u-v|,\text{ for
}b>1. \label{nonlinearity-1}%
\end{equation}
In the case $c=0$, equation (\ref{OKlein}) is called scalar (or non-shifted)
wave equation (for more details see Section \ref{S2}).

The wave and Klein-Gordon equations have been studied by several authors. In
the sequel we review some important results of the literature. Without making
a complete list, we begin by recalling briefly results for equations on the
Euclidean space $\mathbb{R}^{n}$. The existence of global solutions was
analyzed by Georgiev \textit{et al.} \cite{Ge1997}, Ginibre and Velo
\cite{Gi1}, Zhou \cite{Zho1995}, Belchev \textit{et al.} \cite{Bel}, among
others. The nonexistence of global solutions was studied by Sideris \cite{Si}.
Time decays of solutions for wave equations were proved in \cite{Gi2,G3,Str1,
Fec1982}. In particular, Fecher \cite{Fec1982,Fec1984} obtained results on
time decays and established the nonlinear small data scattering for wave and
Klein-Gordon equations in $\mathbb{R}^{3}$. The scattering theory was also
studied by using the time decays in the works \cite{G4, Lin1995, Str2}. After
that, Hidano \cite{Hi1998,Hi2000,Hi2001} obtained the small data scattering
and blow-up theory for nonlinear wave equation on $\mathbb{R}^{n}$ with
$n=3,\,4$ by using the integral representation formula of solutions.
Concerning the nonlinear wave equations, Strichartz estimates were obtained in
\cite{Ge1997,G4,Tao}. Moreover, we quote the results on blow-up \cite{Gl,Jo}
and the life-span of solutions \cite{LiYu,LiZho}.

The well-posedness and scattering theory for wave equations are obtained in
Sobolev spaces by first proving dispersive and Strichartz type estimates in
suitable norms. By following this general spirit, Metcalfe and Taylor
considered wave equations on the three dimensional hyperbolic space
$\mathbb{H}^{3}$ in \cite{MeTa2011,MeTa2012}. They proved dispersive and
Strichartz estimates for the associated linear problem and then obtained the
global well-posedness of solutions in $C^{\infty}(\mathbb{R}\times
\mathbb{H}^{3})$ with smooth compactly supported initial-data $u_{0}%
,\,u_{1}\in C_{0}^{\infty}(\mathbb{H}^{3})$. For general dimension $n\geq2$,
Anker \textit{et al.} \cite{Anker2012,Anker2014} have succeeded to establish
dispersive and Strichartz estimates and obtained the locally well-posed in
$C([-T,T];\,H^{s,\tau}(\mathbb{H}^{n}))\cap C^{1}([-T,T];\,H^{s,\tau
-1}(\mathbb{H}^{n}))$ in \cite{Anker2012} and the globally well-posed in
$C(\mathbb{R};\,H^{\sigma}(\mathbb{H}^{n}))\cap C^{1}(\mathbb{R}%
;\,H^{\sigma-1}(\mathbb{H}^{n}))\cap L^{p}(\mathbb{R};\,L^{q}(\mathbb{H}%
^{n}))$ in \cite{Anker2014}. By employing dispersive and Strichartz estimates
obtained by Anker\textit{ et al.}, French \cite{Fr2012} studied the
scattering for (\ref{OKlein}) with the scattering data space $H^{\sigma
,2}\oplus H^{\sigma-1,2}$. By the same process, Zhang \cite{Zha2015}
established the global well-posedness and scattering theory for wave equations
on nontrapping asymptotically conic manifolds. Tataru \cite{Ta2001} improved
the dispersive estimates for the wave operator in $\mathbb{H}^{n}%
\,\,(n\geq2)$ and use these estimates together with the conformal wave
equations to obtain the global well-posedness of wave equations in
$\mathbb{R}^{n+1}$ with the initial data $(u_{0},u_{1})\in H^{1}\times L^{2}$.
There are some related works about the Strauss conjecture and blow-up for wave euations on hyperbolic spaces \cite{Wang19,Wang2019} and asymptotically flat spacetimes \cite{Wang2014, Me2017} of C. Wang and collaborations.

On the other hand, the existence and scattering theory for the wave equations
on hyperbolic spaces, noncompact manifolds and asymptotically Euclidean
manifolds have been developed by Phillips \textit{et al.} \cite{Phi1987},
Debi\`{e}vre \textit{et al.} \cite{De1992} and Bony and H\"{a}fner
\cite{Bo2008}, respectively, by means of the translation representation for
the unperturbed system and Mourre theory, i.e, the spectral method. The
scattering data spaces in \cite{Phi1987,De1992} are Hilbert spaces as the
completions of $C_{0}^{\infty}(\mathbb{H}^{n})\times C_{0}^{\infty}%
(\mathbb{H}^{n})$ under the energy norms. For the original reference about
this method, we refer the readers to the book of Lax and Phillips
\cite{Lax1967}. Moreover, the geometric scattering and inverse scattering for
the wave equations on some manifolds such as asymptotically Euclidean and
hyperbolic manifolds have been constructed by using the radiation fields (see
for instance \cite{Fri2001, Sa2005, Sa2008} and references therein). In
related works, the scattering theory for nonlinear conformal wave equations on
global hyperbolic space-time was studied by Baez \textit{et al.} \cite{Bae}
and results on the hyperbolic spaces have been studied in
\cite{An2000,Ba2008,Ba2009} for $p$-forms and Schr\"{o}dinger equations.

The purpose of this paper is to establish the local and global well-posedness,
asymptotic stability and scattering for equation (\ref{OKlein}) in a framework based on weak-$L^{p}$ spaces of $\mathbb{H}^{n}$. Moreover, we construct wave operators in such singular setting. The idea of study the well-posedness and scattering in weak-$L^{p}$ spaces were initially developed by Cazenave \textit{et al.} \cite{Ca2001} for Schr\"{o}dinger equations in the
Euclidean space $\mathbb{R}^{n}$ by considering mixed space-time weak-$L^{p}$,
namely $L^{p,\infty}(\mathbb{R}\times\mathbb{R}^{n})$ with $p=\frac
{(b-1)(N+2)}{2}$, and employing Strichartz-type estimates. Then, still
considering Schr\"{o}dinger equations in $\mathbb{R}^{n},$ but employing
dispersive-type estimates, the global well-posedness and asymptotic behavior
of solutions were obtained by Ferreira \textit{et al.} \cite{Fe2009} in a
framework of time polynomial weighted spaces based on the $L^{(p,\infty
)}(\mathbb{R}^{n})$ with $p=b+1$, extending some results obtained in the
$L^{p}$-setting by Cazenave and Weissler \cite{Ca1998, Ca2000}. In this
direction, we also have well-posedness and scattering results for Boussinesq
equations \cite{Fe2011} and wave equations \cite{Liu2009, Fe2017}. In these
works the authors used $L^{(p,r)}$-$L^{(p^{\prime},r)}$-dispersive estimates
on $\mathbb{R}^{n},$ where $1/p+1/p^{\prime}=1$, $1\leq r\leq\infty,$ and
$L^{(p,r)}$ stands for the so-called Lorentz space. The weak-$L^{p}$ space
corresponds to the case $r=\infty.$ The global well-posedness and scattering
in $L^{(p,\infty)}(\mathbb{R}^{n})$ via that approach require the use of
suitable Kato-type classes. For example, in the case of the Schr\"{o}dinger
equation in $\mathbb{R}^{n}$, it is considered the time polynomial weighted
space
\[
\mathcal{H}_{\alpha}=\left\{  u\text{ measurable};\text{ }\sup_{t\in
\mathbb{R}}\text{ }|t|^{\alpha}\left\Vert u(t,\cdot)\right\Vert _{L^{(p,\infty
)}}<\infty\right\}
\]
as well as the initial-data class $\left\{  u_{0}\in\mathcal{S}^{\prime
};\text{ }\sup_{t\in\mathbb{R}}\text{ }|t|^{\alpha}\left\Vert S(t)u_{0}%
\right\Vert _{L^{(p,\infty)}}<\infty\right\}  $, where $S(t)$ is the
Schr\"{o}dinger operator and the power $\alpha$ depends on $p$ via a suitable
relation. We notice that the Ricci curvature of $\mathbb{H}^{n}$ is
$\mathrm{Ric}(\cdot,\cdot)=-(n-1)g(\cdot,\cdot)$, then the dispersive
estimates can be improved in comparison with the ones in Euclidean space
$\mathbb{R}^{n}$. This fact is mentioned in the works by Tataru \cite{Ta2001},
Metcalfe and Taylor \cite{MeTa2011,MeTa2012} and Anker \textit{et al.}
\cite{Anker2012,Anker2014} where the authors provide dispersive estimates with
better decays than the ones in $\mathbb{R}^{n}$.

In this paper we employ the $L^{p}$-dispersive estimates obtained by Tataru
\cite{Ta2001} in order to establish $L^{(p,r)}$-$L^{(p^{\prime},r)}%
$-dispersive estimates in the hyperbolic space $\mathbb{H}^{n}$, where the
decays are of exponential type for large $\left\vert t\right\vert $ and
polynomial for small $\left\vert t\right\vert .$ Using these estimates, we
obtain a global well-posedness result in the mixed time
weighted space $E_{\alpha,\widetilde{\alpha}}$, namely the set of all Bochner
measurable $u:(-\infty,\infty)\rightarrow L^{(b+1,\infty)}$ such that
\begin{equation}
\sup_{\left\vert t\right\vert \geq t_{0}}e^{\alpha\left\vert t\right\vert
}\left\Vert u(\cdot,t)\right\Vert _{(b+1,\infty)}+\sup_{t\in(-t_{0},t_{0}%
)}\left\vert t\right\vert ^{\widetilde{\alpha}}\left\Vert u(\cdot
,t)\right\Vert _{(b+1,\infty)}<\infty, \label{aux-weight-space-1}%
\end{equation}
where $\alpha,\widetilde{\alpha}$ are suitable positive parameters (see
Section \ref{S4} and Theorem \ref{Global}). In view of the exponential weight
in (\ref{aux-weight-space-1}), we show an exponential asymptotic stability result for the global solutions (see Theorem \ref{Stable}), as well as an exponential scattering behavior (see Theorem \ref{SScattering}). Finally, we give the construction of wave operator in Theorem \ref{WaveOp}.

The paper is organized as follows. Section \ref{S2} provides the Klein-Gordon
equations on $\mathbb{H}^{n}$ and the dispersive estimates in the $L^{p}%
$-setting. Section \ref{S3} is devoted to the definition of Lorentz spaces
$L^{(p,q)}$ on $\mathbb{H}^{n}$ and the dispersive estimates of the wave
operator in $L^{(p,q)}$. In Section \ref{S4} we state and prove the main results.

\section{Wave and Klein-Gordon equations on hyperbolic space}

\label{S2} This section is devoted to recall some facts about wave equations on hyperbolic spaces. For further details, the reader is referred to \cite{Anker2012,MeTa2011,Ta2001}.

Let $\mathbb{H}^{n}=\mathbb{H}^{n}(\mathbb{R})$ stand for a real
hyperbolic manifold, where $n\geq2$ is the dimension, endowed with a
Riemannian metric $g$. This space is realized via a hyperboloid in
$\mathbb{R}^{n+1}$ by considering the upper sheet
\[
\left\{  (x_{0},x_{1},...,x_{n})\in\mathbb{R}^{n+1};\text{ }\,x_{0}\geq1\text{
and }x_{0}^{2}-x_{1}^{2}-x_{2}^{2}...-x_{n}^{2}=1\,\right\}
\]
where the metric is given by $dg=-dx_{0}^{2}+dx_{1}^{2}+...+dx_{n}^{2}.$

In geodesic polar coordinates, the hyperbolic manifold $(\mathbb{H}^{n},g)$
can be described as
\[
\mathbb{H}^{n}=\left\{  (\cosh\tau,\omega\sinh\tau),\,\tau\geq0,\omega
\in\mathbb{S}^{n-1}\right\}
\]
with $dg=d\tau^{2}+(\sinh\tau)^{2}d\omega^{2},$ where $d\omega^{2}$ is the
canonical metric on the sphere $\mathbb{S}^{n-1}$. In these coordinates, the
Laplace-Beltrami operator $\Delta_{\mathbb{H}^{n}}$ on $\mathbb{H}^{n}$ can be
expressed as
\[
\Delta_{x}:=\Delta_{\mathbb{H}^{n}}=\partial_{r}^{2}+(n-1)\coth r\partial
_{r}+\sinh^{-2}r\Delta_{\mathbb{S}^{n-1}}.
\]
It is well known that the spectrum of $-\Delta_{x}$ is the half-line
$[\rho^{2},\infty)$, where $\rho=\dfrac{n-1}{2}$.

We consider the nonlinear Klein-Gordon equation on $\mathbb{H}^{n}$
\begin{equation}%
\begin{cases}
\partial_{t}^{2}u(t,x)-\Delta_{x}u(t,x)+cu(t,x)=F(u(t,x)),\\
u(0,x)=u_{0}(x),\,\,\partial_{t}u(t,x)=u_{1}(x),
\end{cases}
\label{Klein}%
\end{equation}
where $c\geq-\dfrac{(n-1)^{2}}{4}$ and the nonlinearity $F(u(t,x))$
satisfies
\begin{equation}
F(0)=0\text{ and }|F(u)-F(v)|\leq C(|u|^{b-1}+|v|^{b-1})|u-v|
\label{Klein-Cond-F}%
\end{equation}
with $b>1$. In particular, we have the two basic cases depending on the value
of $c$:

\begin{itemize}
\item[$\bullet$] If $c=0$, then equation (\ref{Klein}) is the non-shifted wave equation.

\item[$\bullet$] If $c=-\rho^{2}=-\dfrac{(n-1)^{2}}{4}$, then (\ref{Klein}) is
the shifted wave equation.
\end{itemize}

Setting $D=\sqrt{-\Delta_{x}+c}$, Cauchy problem (\ref{Klein}) takes the form
\begin{equation}%
\begin{cases}
\partial_{t}^{2}u(t,x)+D_{x}^{2}u(t,x)=F(u(t,x)),\\
u(0,x)=u_{0}(x),\,\,\partial_{t}u(t,x)=u_{1}(x).
\end{cases}
\label{Klein1}%
\end{equation}
In view of Duhamel's principle, we can formally written (\ref{Klein1}) as
\begin{equation}
u(t)=\dot{W}(t)u_{0}+W(t)u_{1}+\int_{0}^{t}W(t-s)F(u(s))ds,
\label{intergralEq}%
\end{equation}
where $W(t)$ is the wave group. Recall that
\begin{equation}
W(t)=\dfrac{\sin(tD)}{D}\text{ \ and \ }\dot{W}(t)=\cos(tD).
\label{aux-wave-rel-1}%
\end{equation}

The dispersive and Strichartz estimates for the solutions of wave and
Klein-Gordon equations on the hyperbolic manifolds were studied in many works
such as \cite{Anker2012,Anker2014,MeTa2011,MeTa2012,Ta2001} and references
therein. Here, we recall the $L^{p}-L^{p^{\prime}}$-dispersive estimates
obtained by Tataru \cite{Ta2001}.

\begin{proposition}
\label{Dispersive} Let $n\geq2$, $2\leq p\leq\frac{2(n+1)}{n-1}$ and $\frac
{1}{p}+\frac{1}{p^{\prime}}=1$. Then, there exists a constant $C>0$
(independent of $t$) such that
\begin{equation}
\left\Vert W(t)g\right\Vert _{L^{p}}\leq C\phi_{p}(t)\left\Vert g\right\Vert
_{L^{p^{\prime}}}, \label{group1}%
\end{equation}
for all $g\in L^{p^{\prime}}(\mathbb{H}^{n}),$ where $\phi_{p}(t)=\frac
{(1+\left\vert t\right\vert )^{\frac{2}{p}}}{(\sinh\left\vert t\right\vert
)^{\frac{n-1}{2}\left(  1-\frac{2}{p}\right)  }}$.
\end{proposition}

\begin{remark}
\label{Re}

\begin{itemize}
\item[$(i)$] We notice that the above proposition is a consequence of
\cite[Theorem 3]{Ta2001} by considering $2s=\frac{n+1}{2}\left(  1-\frac{2}%
{p}\right)  \leq1$. This condition is equivalent to $p\leq\frac{2(n+1)}{n-1}$.

\item[$(ii)$] Proceeding as in the proof of Theorem 3 in \cite{Ta2001}, we
also have the estimate
\begin{equation}
\left\Vert \frac{\dot{W}(t)}{D}g\right\Vert _{L^{p}}=\left\Vert \frac
{\cos(tD)}{D}g\right\Vert _{L^{p}}\leq C\phi_{p}(t)\left\Vert g\right\Vert
_{L^{p^{\prime}}}, \label{Est-disp-2}%
\end{equation}
for all $g\in L^{p^{\prime}}(\mathbb{H}^{n})$.
\end{itemize}
\end{remark}

\section{Lorentz spaces and interpolation estimates on hyperbolic spaces}

\label{S3}

Let $\Omega$ be a subset of $\mathbb{H}^{n}$. For $0<p<\infty,$ denote by
$L^{p}(\Omega)$ the space of all $L^{p}$-integrable functions on $\Omega.$ Let
$1<p_{1}<p_{2}\leq\infty,$ $\theta\in(0,1)$ and $1\leq
r\leq\infty$ with $\dfrac{1}{p}=\dfrac{1-\theta}{p_{1}}+\dfrac{\theta
}{p_{2}}.$ The Lorentz space $L^{(p,r)}(\Omega)$ is defined as the
interpolation space $\left(  L^{p_{1}},L^{p_{2}}\right)  _{\theta,r}=L^{(p,r)}
$ with the natural norm $\left\Vert \cdot\right\Vert _{(p,r)}$ induced by the
functor $(\cdot,\cdot)_{\theta,r}.$ In particular, $L^{p}(\Omega
)=L^{(p,p)}(\Omega)$ and $L^{(p,\infty)}$ is the so-called weak-$L^{p}$ space
or the Marcinkiewicz space on $\Omega$.

For $1\leq q_{1}\leq p\leq q_{2}\leq\infty,$ we have the
following relation
\begin{equation}
L^{(p,1)}(\Omega)\subset L^{(p,q_{1})}(\Omega)\subset L^{p}(\Omega)\subset
L^{(p,q_{2})}(\Omega)\subset L^{(p,\infty)}(\Omega). \label{Inclusion}%
\end{equation}
Let $1<p_{1},p_{2},p_{3}\leq\infty$ and $1\leq r_{1},r_{2}%
,r_{3}\leq\infty$ be such that $\dfrac{1}{p_{3}}=\dfrac{1}{p_{1}}%
+\dfrac{1}{p_{2}}$ and $\dfrac{1}{r_{1}}+\dfrac{1}{r_{2}}\geq\dfrac
{1}{r_{3}}$. We have the H\"{o}lder inequality%
\begin{equation}
\left\Vert fg\right\Vert _{{(p_{3},r_{3})}}\leq C\left\Vert f\right\Vert
_{{(p_{1},r_{1})}}\left\Vert g\right\Vert _{{(p_{2},r_{2})}}, \label{Holder}%
\end{equation}
where $C>0$ is a constant independent of $f$ and $g$. Moreover, for
$1<p_{1}<p_{2}\leq\infty$, $0<\theta<1$, $1\leq r_{1},r_{2}%
,r\leq\infty,$ and $\dfrac{1}{p}=\dfrac{1-\theta}{p_{1}}+\dfrac{\theta
}{p_{2}},$ by reiteration theorem (see \cite[Theorem 3.5.3 ]{BeLo}), we have
the interpolation property
\begin{equation}
\left(  L^{(p_{1},r_{1})},L^{(p_{2},r_{2})}\right)  _{\theta,r}=L^{(p,r)}.
\label{interp1}%
\end{equation}

Now we extend the dispersive estimates in Proposition \ref{Dispersive} to the
framework of Lorentz spaces.

\begin{lemma}
\label{DisperInter} Let $n\geq2$, $1\leq r\leq\infty$, $2<p<\frac{2(n+1)}%
{n-1}$ and $\frac{1}{p}+\frac{1}{p^{\prime}}=1$. Then, there exists a constant
$C>0$ (independent of $t$) such that
\begin{equation}
\left\Vert W(t)g\right\Vert _{{(p,r)}}+\,\left\Vert \frac{\dot{W}(t)}%
{D}g\right\Vert _{{(p,r)}}\leq C\phi_{p}(t)\left\Vert g\right\Vert
_{{(p^{\prime},r)}}, \label{group2}%
\end{equation}
for all $g\in L^{p^{\prime},r}(\mathbb{H}^{n}),$ where $\phi_{p}%
(t)=\frac{(1+\left\vert t\right\vert )^{\frac{2}{p}}}{(\sinh\left\vert
t\right\vert )^{\frac{n-1}{2}\left(  1-\frac{2}{p}\right)  }}$.
\end{lemma}

\textbf{Proof.} The proof follows by means of an interpolation argument. For
that, let $p_{1},p_{2},p\in(2,\frac{2(n+1)}{n-1})$ and $\theta\in(0,1)$ be
such that $p^{-1}=\theta p_{1}^{-1}+(1-\theta)p_{2}^{-1}.$ Employing
(\ref{group1}) and (\ref{Est-disp-2}), and recalling that $L^{p}=L^{(p,p)},$
we have%
\begin{equation}
\left\Vert W(t)g\right\Vert _{(p_{i},p_{i})}+\,\left\Vert \frac{\dot{W}(t)}%
{D}g\right\Vert _{(p_{i},p_{i})}\leq C_{i}\phi_{p_{i}}(t)\left\Vert
g\right\Vert _{(p_{i}^{\prime},p_{i}^{\prime})}, \label{aux-est-interp-1}%
\end{equation}
where $C_{i}>0$ is independent of $t.$ Recalling that the interpolation
functor $\left(  \cdot,\cdot\right)  _{\theta,r}$ is of exponent $\theta$ and
$\left(  L^{p_{1}},L^{p_{2}}\right)  _{\theta,r}=L^{(p,r)}$ (see \cite{BeLo}),
estimate (\ref{aux-est-interp-1}) leads us to
\begin{equation}
\left\Vert W(t)g\right\Vert _{(p,r)}+\,\left\Vert \frac{\dot{W}(t)}%
{D}g\right\Vert _{(p,r)}\leq K(C_{1}\phi_{p_{1}}(t))^{\theta}(C_{2}\phi
_{p_{2}}(t))^{1-\theta}\left\Vert g\right\Vert _{(p^{\prime},r)},
\label{aux-est-interp-2}%
\end{equation}
where $K>0$ is independent of $t$. Next, taking $C=K\left(  C_{1}\right)
^{\theta}(C_{2})^{1-\theta}$ and using $p^{-1}=\theta p_{1}^{-1}%
+(1-\theta)p_{2}^{-1},$ note that
\begin{align*}
K(C_{1}\phi_{p_{1}}(t))^{\theta}(C_{2}\phi_{p_{2}}(t))^{1-\theta}  &
=K\left(  C_{1}\right)  ^{\theta}(C_{2})^{1-\theta}\frac{(1+\left\vert
t\right\vert )^{\frac{2\theta}{p_{1}}}}{(\sinh\left\vert t\right\vert
)^{\frac{n-1}{2}\theta\left(  1-\frac{2}{p_{1}}\right)  }}\frac{(1+\left\vert
t\right\vert )^{\frac{2(1-\theta)}{p_{2}}}}{(\sinh\left\vert t\right\vert
)^{\frac{n-1}{2}(1-\theta)\left(  1-\frac{2}{p_{2}}\right)  }}\\
&  =C\frac{(1+\left\vert t\right\vert )^{\frac{2}{p}}}{(\sinh\left\vert
t\right\vert )^{\frac{n-1}{2}\left(  1-\frac{2}{p}\right)  }},
\end{align*}
which, together with (\ref{aux-est-interp-2}), gives the desired estimate. \fin

\section{Well-posedness, stability and scattering}

\label{S4}

\noindent In order to establish the well-posedness, we define some suitable functional spaces as follows. Let $T>0$, $1<b<\infty$,
$\eta\geq0$, and denote by $E_{\eta}^{T}$ the set of all Bochner measurable
$u:(-T,T)\rightarrow L^{(b+1,\infty)}(\mathbb{H}^{n})$ such that
\[
\left\Vert u\right\Vert _{E_{\eta}^{T}}=\sup_{t\in(-T,T)}|t|^{\eta
}\left\Vert u\right\Vert _{(b+1,\infty)(\mathbb{H}^n)}<\infty.
\]
The pair $(E_{\eta}^{T},\left\Vert \cdot\right\Vert _{E_{\eta}^{T}})$ is a
Banach space.

Now, for a fixed $t_{0}\geq1$ (see (\ref{T0})), $\alpha,\widetilde{\alpha}>0$, $1<b<\infty$, and $1\leq d\leq\infty$, consider the
space $E_{\alpha,\widetilde{\alpha}}^{d}$ of all Bochner measurable
$u:(-\infty,\infty)\rightarrow L^{(b+1,d)}$ satisfying
\begin{equation}
\left\Vert u\right\Vert _{E_{\alpha,\widetilde{\alpha}}^{d}}=\sup_{\left\vert
t\right\vert \geq t_{0}}e^{\alpha\left\vert t\right\vert }\left\Vert
u(t)\right\Vert _{(b+1,d)}+\sup_{t\in(-t_{0},t_{0})}\left\vert t\right\vert
^{\widetilde{\alpha}}\left\Vert u(t)\right\Vert _{(b+1,d)}<\infty
.\label{space-1}%
\end{equation}
The space $E_{\alpha,\widetilde{\alpha}}^{d}$ endowed with $\left\Vert
\cdot\right\Vert _{E_{\alpha,\widetilde{\alpha}}^{d}}$ is a Banach space.
Also, we consider the initial-data space $\mathcal{I}_{0}^{d}$ as the set of
all pairs $(u_{0},u_{1})\in\mathcal{S}^{\prime
}(\mathbb{H}^{n})\times\mathcal{S}^{\prime}(\mathbb{H}^{n})$ such that
\begin{equation}
\left\Vert (u_{0},u_{1})\right\Vert _{\mathcal{I}_{0}^{d}}=\left\Vert \dot
{W}(t)u_{0}+W(t)u_{1}\right\Vert _{E_{\alpha,\widetilde{\alpha}}^{d}}%
<\infty.\label{initial}%
\end{equation}
In the case $d=\infty,$ we use the notation%
\begin{equation}
E_{\alpha,\widetilde{\alpha}}=E_{\alpha,\widetilde{\alpha}}^{\infty}\text{
\ \ \ and \ }\mathcal{I}_{0}=\mathcal{I}_{0}^{\infty}.\label{space-infty}%
\end{equation}

\subsection{Global well-posedness}

\begin{theorem}
\label{Global}(Global-in-time solution) \noindent

\begin{itemize}
\item[$(i)$] (Well-posedness). Let $\beta=\frac{n-1}{2}\left(  1-\frac{2}{b+1}\right)$, $0<\sigma<\beta$, and $\frac{n+1+\sigma+ \sqrt{(n+1+\sigma)^{2}+8(n-1-\sigma)}}{2(n-1-\sigma)} < b < \frac{n+3}{n-1}$ satisfy $1- (\beta-\sigma) = (b-1)\alpha$ and $1-\beta = (b-1)\widetilde{\alpha}$.
If $(u_{0},u_{1})\in\mathcal{I}_{0}$ with $\left\Vert \left(  u_{0}%
,u_{1}\right)  \right\Vert _{\mathcal{I}_{0}}\leq\varepsilon$ for some
$\varepsilon>0$ small enough, then equation \eqref{intergralEq} has a unique
global-in-time mild solution $u\in E_{\alpha,\widetilde{\alpha}}$ such that
$\left\Vert u\right\Vert _{E_{\alpha,\widetilde{\alpha}}}\leq2\varepsilon$.

\item[$(ii)$] ($L^{p,q}$-regularity). Let $1\leq d\leq\infty$ and $0\leq
h < 1-b\alpha$. Suppose further that%
\begin{equation}
\Gamma_{1,h}^{d}:=\sup_{\left\vert t\right\vert \geq t_{0}}e^{(\alpha
+h)\left\vert t\right\vert }\left\Vert \dot{W}(t)u_{0}+W(t)u_{1}\right\Vert
_{(b+1,d)}+\sup_{t\in(-t_{0},t_{0})}\left\vert t\right\vert ^{\widetilde
{\alpha}+h}\left\Vert \dot{W}(t)u_{0}+W(t)u_{1}\right\Vert _{(b+1,d)}%
<\infty.\label{cond-d-h}%
\end{equation}
Then, there exists $\varepsilon>0$ such that the previous solution satisfies
\begin{equation}
\sup_{\left\vert t\right\vert \geq t_{0}}e^{(\alpha+h)\left\vert t\right\vert
}\left\Vert u(t)\right\Vert _{(b+1,d)}+\sup_{t\in(-t_{0},t_{0})}\left\vert
t\right\vert ^{\widetilde{\alpha}+h}\left\Vert u(t)\right\Vert _{(b+1,d)}%
<\infty,\label{Globalregular}%
\end{equation}
provided that $\left\Vert \left(  u_{0},u_{1}\right)  \right\Vert
_{\mathcal{I}_{0}}\leq\varepsilon$.
\end{itemize}
\end{theorem}

\begin{remark}
\label{Local}(Local-in-time solution). \noindent\ With an adaptation on the
arguments, we can obtain a local version of the above well-posedness result,
regardless of the initial-data size. More precisely, let $1<b<\frac{n+1+\sqrt{n^{2}+10n-7}}{2(n-1)}$ and $0<\beta=\frac{n-1}{2}\left(  1-\frac{2}{b+1}\right)  <1$. If $(Du_{0},u_{1})\in
L^{(\frac{b+1}{b},\infty)}(\mathbb{H}^n)\times L^{(\frac{b+1}{b},\infty)}(\mathbb{H}^n)$, then there
exists $0<T<\infty$, such that (\ref{intergralEq}) has a unique solution $u\in
E_{\eta}^{T}$ with $\eta=\beta$. Moreover, if $(Du_{0},u_{1})\in L^{(\frac{b+1}{b},d)}(\mathbb{H}^n)\times
L^{(\frac{b+1}{b},d)}(\mathbb{H}^n)$ with $1\leq d\leq\infty$, then the solution $u$
satisfies
\begin{equation}
\sup_{t\in(-T,T)}|t|^{\beta}\left\Vert u\right\Vert _{(b+1,d)}<\infty
.\label{Localregular}%
\end{equation}

\end{remark}

\textbf{Proof of Theorem \ref{Global}.}

\textbf{Item (i).} First note that we can choose $t_{0}\geq1$ such that the
function $\phi_{p}(t)$ in (\ref{group2}) satisfies
\begin{equation}
\phi_{p}(t)\leq%
\begin{cases}
\displaystyle C\left\vert t\right\vert ^{\frac{2}{p}}e^{-\beta_{p}\left\vert
t\right\vert }, & \text{for \ }\left\vert t\right\vert \geq t_{0},\\
C\left\vert t\right\vert ^{-\beta_{p}}, & \text{for \ }0<\left\vert t\right\vert
<t_{0},
\end{cases}
\label{T0}%
\end{equation}
where $\beta_{p}=\frac{n-1}{2}\left(1-\frac{2}{p}\right)$. Also, denote
\begin{equation}
\mathcal{T}(u)=\int_{0}^{t}W(t-s)F(u(s))ds. \label{Op1}%
\end{equation}

From the conditions on $\sigma,\, b,\, \alpha$ and $\widetilde{\alpha}$, we have the relations
\begin{align*}
&2<b+1<\frac{2(n+1)}{n-1},\, \widetilde{\alpha} = \frac{n-1}{b^2-1}-\frac{n-3}{2(b-1)},\, \alpha = \widetilde{\alpha} + \frac{\sigma}{b-1},\\
&\text{ \ } \\
&0<\beta<1,\, 0<b\widetilde{\alpha}<b\alpha<1,\; \hbox{  and  } \; 0<1-b\alpha=(\beta-\sigma)-\alpha.
\end{align*}

In view of the time symmetry of the wave group (\ref{Op1}) and estimates
(\ref{group2}) and (\ref{T0}), we can assume $t>0$ without loss of generality.
Moreover, we consider three cases for the time variable $t$ when estimating
the operator $\mathcal{T}$.

\underline{$\bullet$ For $t_{0}-\delta\leq t\leq t_{0}+\delta$, where $0<\delta<t_0$:} In this case, letting $p=b+1$ in (\ref{T0}), we have $\beta_{p}=\beta$, $e^{\alpha t}\sim t^{\widetilde{\alpha}}$, and $t^{\frac{2}{b+1}}e^{-\beta t}\sim t^{-\beta}.$ Then, it is not difficult to see that
\begin{equation}
\sup_{t\in(t_{0}-\delta,t_{0}+\delta)}\left(  \max\{e^{\alpha\left\vert t\right\vert
},\left\vert t\right\vert ^{\widetilde{\alpha}}\}\left\Vert \mathcal{T}%
(u)-\mathcal{T}(v)\right\Vert _{(b+1,\infty)}\right)  \leq C_{1}\left\Vert
u-v\right\Vert _{E_{\alpha,\widetilde{\alpha}}}(\left\Vert u\right\Vert
_{E_{\alpha,\widetilde{\alpha}}}^{b-1}+\left\Vert v\right\Vert _{E_{\alpha
,\widetilde{\alpha}}}^{b-1}). \label{CCCore}%
\end{equation}

\underline{$\bullet$ For $t\geq t_{0}+\delta$:} First note that
$1-\beta-b\widetilde{\alpha}=-\widetilde{\alpha},$ $b\widetilde{\alpha}<1$ and
$b\alpha<1.$ Then, using (\ref{T0}), Lemma \ref{DisperInter}, Remark \ref{Re}
(ii) and H\"{o}lder's inequality, we can estimate
\begin{align}
\left\Vert \mathcal{T}(u)-\mathcal{T}(v)\right\Vert _{(b+1,\infty)} &
\leq\int_{0}^{t}\left\Vert W(t-s)[F(u)-F(v)](s)\right\Vert _{(b+1,\infty
)}ds \nonumber\\
&  =\int_{0}^{t}\left\Vert W(s)[F(u)-F(v)](t-s)\right\Vert _{(b+1,\infty
)}ds \nonumber\\
&  \leq C\int_{0}^{t}\phi_{b+1}(s)\left\Vert [F(u)-F(v)](t-s)\right\Vert
_{(\frac{b+1}{b},\infty)}ds \nonumber\\
&  \leq C\int_{0}^{t}\phi_{b+1}(s)\left\Vert (u-v)(t-s)\right\Vert
_{(b+1,\infty)}\left(  \left\Vert u(t-s)\right\Vert _{(b+1,\infty)}%
^{b-1}+\left\Vert v(t-s)\right\Vert _{(b+1,\infty)}^{b-1}\right)
ds \nonumber\\
&  \leq C\int_{t_{0}}^{t}s^{\frac{2}{b+1}}e^{-\frac{n-1}{2}\left(
1-\frac{2}{b+1}\right)  s}\left\Vert (u-v)(t-s)\right\Vert _{(b+1,\infty
)}\left(  \left\Vert u(t-s)\right\Vert _{(b+1,\infty)}^{b-1}+\left\Vert
v(t-s)\right\Vert _{(b+1,\infty)}^{b-1}\right)  ds \nonumber\\
&  +C\int_{0}^{t_{0}}s^{-\frac{n-1}{2}\left(  1-\frac{2}{b+1}\right)
}\left\Vert (u-v)(t-s)\right\Vert _{(b+1,\infty)}\left(  \left\Vert
u(t-s)\right\Vert _{(b+1,\infty)}^{b-1}+\left\Vert v(t-s)\right\Vert
_{(b+1,\infty)}^{b-1}\right)  ds \nonumber\\
&  =I_{1}+I_{2}.\label{ine}%
\end{align}
Considering $p=b+1$ and $\beta_{p}=\beta$ in (\ref{T0}), and taking $\delta=\dfrac{t_{0}}{2}$, we estimate the integral $I_{1}$ as follows
\begin{align}
I_{1} &  =C\int_{t_{0}}^{t}s^{\frac{2}{b+1}}e^{-\frac{n-1}{2}\left(
1-\frac{2}{b+1}\right)  s}\left\Vert (u-v)(t-s)\right\Vert _{(b+1,\infty
)}\left(  \left\Vert u(t-s)\right\Vert _{(b+1,\infty)}^{b-1}+\left\Vert
v(t-s)\right\Vert _{(b+1,\infty)}^{b-1}\right)  ds \nonumber\\
&  =C\int_{t_{0}}^{t-\delta}s^{\frac{2}{b+1}}e^{-\frac{n-1}{2}\left(
1-\frac{2}{b+1}\right)  s}\left\Vert (u-v)(t-s)\right\Vert _{(b+1,\infty
)}\left(  \left\Vert u(t-s)\right\Vert _{(b+1,\infty)}^{b-1}+\left\Vert
v(t-s)\right\Vert _{(b+1,\infty)}^{b-1}\right)  ds \nonumber\\
&  +C\int_{t-\delta}^{t}s^{\frac{2}{b+1}}e^{-\frac{n-1}{2}\left(  1-\frac
{2}{b+1}\right)  s}\left\Vert (u-v)(t-s)\right\Vert _{(b+1,\infty)}\left(
\left\Vert u(t-s)\right\Vert _{(b+1,\infty)}^{b-1}+\left\Vert
v(t-s)\right\Vert _{(b+1,\infty)}^{b-1}\right)  ds \nonumber\\
&  \leq C\sup_{t>\delta}e^{\alpha t}\left\Vert (u-v)(t)\right\Vert
_{(b+1,\infty)}\left(  \sup_{t>\delta}e^{(b-1)\alpha t}\left\Vert
u(t)\right\Vert _{(b+1,\infty)}^{b-1}+\sup_{t>\delta}e^{(b-1)\alpha
t}\left\Vert v(t)\right\Vert _{(b+1,\infty)}^{b-1}\right) \nonumber\\
&  \times\int_{t_{0}}^{t-\delta}s^{\frac{2}{b+1}}e^{-\beta s}e^{-b\alpha
(t-s)}ds \nonumber\\
&  +C\sup_{0<t<\delta}t^{\widetilde{\alpha}}\left\Vert (u-v)(t)\right\Vert
_{(b+1,\infty)}\left(  \sup_{0<t<\delta}t^{(b-1)\widetilde{\alpha}}\left\Vert
u(t)\right\Vert _{(b+1,\infty)}^{b-1}+\sup_{0<t<\delta}t^{(b-1)\widetilde
{\alpha}}\left\Vert v(t)\right\Vert _{(b+1,\infty)}^{b-1}\right)  \nonumber\\
&  \times\int_{t-\delta}^{t}s^{\frac{2}{b+1}}e^{-\beta s}(t-s)^{-b\widetilde
{\alpha}}ds \nonumber\\
&  \leq{C}e^{-\alpha t}\left\Vert u-v\right\Vert _{E_{\alpha,\widetilde
{\alpha}}}(\left\Vert u\right\Vert _{E_{\alpha,\widetilde{\alpha}}}%
^{b-1}+\left\Vert v\right\Vert _{E_{\alpha,\widetilde{\alpha}}}^{b-1}%
),\label{ine1}%
\end{align}
where we have used (\ref{CCCore}), $1-(\beta-\sigma)-b\alpha=-\alpha$, $0<\sigma <\beta<1$, and
\begin{align*}
\int_{t_{0}}^{t-\delta}s^{\frac{2}{b+1}}e^{-\beta s}e^{-b\alpha(t-s)}ds &
\leq\int_{t_{0}}^{t}s^{\frac{2}{b+1}}e^{-\beta s}e^{-b\alpha(t-s)}ds \leq\int_{t_{0}}^{t} \left( s^{\frac{2}{b+1}}e^{-\sigma s} \right) e^{-(b-\sigma)s}e^{-b\alpha(t-s)}ds\\
&  \leq C\int_{t_{0}}^{t}e^{-(\beta-\sigma)s}e^{-b\alpha(t-s)}ds \hbox{   (because  } s^{\frac{2}{b+1}}e^{-\sigma s}\leq C)\\
&  \leq C\int_{t_{0}}^{t}e^{(b\alpha-\alpha-1)s}e^{-b\alpha(t-s)}ds\\
&  \leq Ce^{-b\alpha t}\int_{t_{0}}^{t}e^{(2b\alpha-\alpha-1)s}ds\\
&  \leq\left\vert \frac{Ce^{-b\alpha t}}{2b\alpha-\alpha-1}\left(
e^{(2b\alpha-\alpha-1)t}-e^{(2b\alpha-\alpha-1)t_{0}}\right)  \right\vert \\
&  \leq Ce^{-\alpha t} \hbox{   (because   } b\alpha<1)
\end{align*}
and%
\begin{align*}
\int_{t-\delta}^{t}s^{\frac{2}{b+1}}e^{-\beta s}(t-s)^{-b\widetilde{\alpha}}ds
&  \leq\int_{t-\delta}^{t}s^{\frac{2}{b+1}}e^{-\beta s}(t-s)^{-b\widetilde
{\alpha}}ds\\
&  \leq C\int_{t-\delta}^{t}e^{-(\beta-\sigma)s}(t-s)^{-b\widetilde{\alpha}%
}ds\\
&  \leq C\int_{t-\delta}^{t}e^{(b\alpha-\alpha-1)s}(t-s)^{-b\widetilde{\alpha
}}ds\\
&  =C\int_{0}^{\delta}e^{(b\alpha-\alpha-1)(t-s)}s^{-b\widetilde{\alpha}}ds\\
&  =C\int_{0}^{\delta/2}e^{(b\alpha-\alpha-1)(t-s)}s^{-b\widetilde{\alpha}%
}ds+C\int_{\delta/2}^{\delta}e^{(b\alpha-\alpha-1)(t-s)}s^{-b\widetilde
{\alpha}}ds\\
&  \leq Ce^{(b\alpha-\alpha-1)\left(  t-\frac{\delta}{2}\right)  }\int
_{0}^{\delta/2}s^{-b\widetilde{\alpha}}ds+C(\frac{\delta}{2})^{-b\widetilde
{\alpha}}\int_{\delta/2}^{\delta}e^{(b\alpha-\alpha-1)(t-s)}ds\\
&  \leq Ce^{-\alpha t} \hbox{   (because   } b\widetilde{\alpha}<b\alpha<1).
\end{align*}
The remaining term $I_{2}$ can be estimated as
\begin{align}
I_{2} &  =C\int_{0}^{t_{0}}s^{-\frac{n-1}{2}\left(  1-\frac{2}{b+1}\right)
}\left\Vert (u-v)(t-s)\right\Vert _{b+1,\infty}\left(  \left\Vert
u(t-s)\right\Vert _{b+1,\infty}^{b-1}+\left\Vert v(t-s)\right\Vert
_{b+1,\infty}^{b-1}\right)  ds \nonumber\\
&  \leq C\sup_{t>\delta}e^{\alpha t}\left\Vert (u-v)(t)\right\Vert
_{(b+1,\infty)}\left(  \sup_{t>\delta}e^{(b-1)\alpha t}\left\Vert
u(t)\right\Vert _{(b+1,\infty)}^{b-1}+\sup_{t>\delta}e^{(b-1)\alpha
t}\left\Vert v(t)\right\Vert _{(b+1,\infty)}^{b-1}\right)  \newline \nonumber\\
&  \times\int_{0}^{t_{0}}s^{-\beta}e^{-b\alpha(t-s)}ds \nonumber\\
&  \leq Ce^{-\alpha t}\left\Vert u-v\right\Vert _{E_{\alpha,\widetilde{\alpha
}}}(\left\Vert u\right\Vert _{E_{\alpha,\widetilde{\alpha}}}^{b-1}+\left\Vert
v\right\Vert _{E_{\alpha,\widetilde{\alpha}}}^{b-1}),\label{ine2}%
\end{align}
where we have used (\ref{CCCore}) and the inequality
\begin{align*}
\int_{0}^{t_{0}}s^{-\beta}e^{-b\alpha(t-s)}ds &  =\int_{0}^{t_{0}/2}s^{-\beta
}e^{-b\alpha(t-s)}ds+\int_{t_{0}/2}^{t_{0}}s^{-\beta}e^{-b\alpha(t-s)}ds\\
&  \leq e^{-b\alpha\left(  t-\frac{t_{0}}{2}\right)  }\frac{1}{1-\beta}%
(\frac{t_{0}}{2})^{1-\beta}+(\frac{t_{0}}{2})^{-\beta}\frac{1}{b\alpha}\left(
e^{-b\alpha(t-t_{0})}-e^{-b\alpha(t-\frac{t_{0}}{2})}\right)  \\
&  \leq Ce^{-\alpha t},
\end{align*}
because $t> t_{0}$, $b>1$ and $0<\beta<1$. Combining (\ref{ine}), (\ref{ine1}) and
(\ref{ine2}), we arrive at
\begin{equation}
\sup_{t\geq t_{0}+\delta}e^{\alpha t}\left\Vert \mathcal{T}(u)-\mathcal{T}%
(v)\right\Vert _{(b+1,\infty)}\leq C_{2}\left\Vert u-v\right\Vert
_{E_{\alpha,\widetilde{\alpha}}}(\left\Vert u\right\Vert _{E_{\alpha
,\widetilde{\alpha}}}^{b-1}+\left\Vert v\right\Vert _{E_{\alpha,\widetilde
{\alpha}}}^{b-1}).\label{Core}%
\end{equation}

\underline{$\bullet$ For $0<t\leq t_{0}$:} Employing (\ref{T0}), Lemma
\ref{DisperInter}, Remark \ref{Re} (ii), and H\"{o}lder's inequality, we
obtain that
\begin{align}
\left\Vert \mathcal{T}(u)-\mathcal{T}(v)\right\Vert _{(b+1,\infty)}  &
\leq\int_{0}^{t}\left\Vert W(t-s)[F(u)-F(v)](s)\right\Vert _{(b+1,\infty
)}ds \nonumber\\
&  \leq C\int_{0}^{t}\phi_{b+1}(s)\left\Vert [F(u)-F(v)](t-s)\right\Vert
_{(\frac{b+1}{b},\infty)}ds \nonumber\\
&  \leq C\int_{0}^{t}\phi_{b+1}(s)\left\Vert (u-v)(t-s)\right\Vert
_{(b+1,\infty)}\left(  \left\Vert u(t-s)\right\Vert _{(b+1,\infty)}%
^{b-1}+\left\Vert v(t-s)\right\Vert _{(b+1,\infty)}^{b-1}\right)
ds \nonumber\\
&  \leq C\int_{0}^{t}s^{-\frac{n-1}{2}\left(  1-\frac{2}{b+1}\right)
}\left\Vert (u-v)(t-s)\right\Vert _{(b+1,\infty)}\left(  \left\Vert
u(t-s)\right\Vert _{(b+1,\infty)}^{b-1}+\left\Vert v(t-s)\right\Vert
_{(b+1,\infty)}^{b-1}\right)  ds \nonumber\\
&  \leq C\int_{0}^{t}s^{-\beta}(t-s)^{-b\widetilde{\alpha}}ds \sup_{t\in(0,t_{0})}t^{\widetilde{\alpha}}\left\Vert
(u-v)(t)\right\Vert _{(b+1,\infty)} \nonumber\\
&  \times\left(\sup_{t\in(0,t_{0})}t^{(b-1)\widetilde{\alpha}}\left\Vert u(t)\right\Vert _{(b+1,\infty)}%
^{b-1}+\sup_{t\in(0,t_{0})}t^{(b-1)\widetilde{\alpha}}\left\Vert
v(t)\right\Vert _{(b+1,\infty)}^{b-1}\right) \nonumber\\
&  \leq C_{3}t^{-\widetilde{\alpha}}\left\Vert u-v\right\Vert _{E_{\alpha
,\widetilde{\alpha}}}(\left\Vert u\right\Vert _{E_{\alpha,\widetilde{\alpha}}%
}^{b-1}+\left\Vert v\right\Vert _{E_{\alpha,\widetilde{\alpha}}}^{b-1}),
\label{1ine}%
\end{align}
where we have used in (\ref{1ine}) that $\int_{0}^{t}s^{-\beta}%
(t-s)^{-b\widetilde{\alpha}}ds=Ct^{-\widetilde{\alpha}}.$ Therefore, we get
\begin{equation}
\sup_{t\in(0,t_{0})}t^{\widetilde{\alpha}}\left\Vert \mathcal{T}%
(u)-\mathcal{T}(v)\right\Vert _{(b+1,\infty)}\leq C_{3}\left\Vert
u-v\right\Vert _{E_{\alpha,\widetilde{\alpha}}}(\left\Vert u\right\Vert
_{E_{\alpha,\widetilde{\alpha}}}^{b-1}+\left\Vert v\right\Vert _{E_{\alpha
,\widetilde{\alpha}}}^{b-1}). \label{CCore}%
\end{equation}

Next, for some $\varepsilon>0$, we are going to show that $\mathcal{T}(u)$ is
a contraction in the closed ball $B(0,2\varepsilon)=\left\{  u\in
E_{\alpha,\widetilde{\alpha}};\left\Vert u\right\Vert _{E_{\alpha
,\widetilde{\alpha}}}\leq2\varepsilon\right\}  \subset E_{\alpha
,\widetilde{\alpha}}$ provided that $\left\Vert (u_{0},u_{1})\right\Vert
_{\mathcal{I}_{0}}\leq\varepsilon$. For that, set $y=\dot{W}%
(t)u_{0}+W(t)u_{1}$ and
\begin{equation}
\Phi(u)(t)=y+\int_{0}^{t}W(t-s)F(u(s))ds.\label{Op2}%
\end{equation}
Using (\ref{CCCore}), (\ref{Core}), (\ref{CCore}) and letting $\mathcal{K}%
=C_{1}+C_{2}+C_{3},$ we have that
\begin{align*}
\left\Vert \Phi(u)\right\Vert _{E_{\alpha,\widetilde{\alpha}}} &
\leq\left\Vert y\right\Vert _{E_{\alpha,\widetilde{\alpha}}}+\int_{0}%
^{t}\left\Vert W(t-s)[F(u(s))]ds\right\Vert _{E_{\alpha,\widetilde{\alpha}}}\\
&  \leq\left\Vert y\right\Vert _{E_{\alpha,\widetilde{\alpha}}%
}+\mathcal{K}\left\Vert u\right\Vert _{E_{\alpha,\widetilde{\alpha}}}^{b}\\
&  \leq\varepsilon+\mathcal{K}2^{b}\varepsilon^{b}\\
&  \leq2\varepsilon
\end{align*}

and%
\begin{align*}
\left\Vert \Phi(u)-\Phi(v)\right\Vert _{E_{\alpha,\widetilde{\alpha}}} &
\leq\int_{0}^{t}\left\Vert W(t-s)[F(u(s))-F(v(s))]ds\right\Vert
_{E_{\alpha,\widetilde{\alpha}}}\\
&  \leq C\left\Vert u-v\right\Vert _{E_{\alpha,\widetilde{\alpha}}%
}(\left\Vert u\right\Vert _{E_{\alpha,\widetilde{\alpha}}}^{b-1}+\left\Vert
v\right\Vert _{E_{\alpha,\widetilde{\alpha}}}^{b-1})\\
&  \leq\mathcal{K}2^{b}\varepsilon^{b-1}\left\Vert u-v\right\Vert
_{E_{\alpha,\widetilde{\alpha}}},
\end{align*}
for all $u,v\in B(0,2\varepsilon),$ provided that $\mathcal{K}2^{b}%
\varepsilon^{b-1}<1$. It follows that $\Phi:B(0,2\varepsilon)\rightarrow
B(0,2\varepsilon)$ is a contraction. Then, in view of the Banach fixed point
theorem, we can conclude the existence of a unique solution $u\in
E_{\alpha,\widetilde{\alpha}}$ for (\ref{intergralEq}) such that $\left\Vert
u\right\Vert _{E_{\alpha,\widetilde{\alpha}}}\leq2\varepsilon.$

\textbf{Item (ii).} Let $0\leq h<1-b\alpha,$ $1\leq d\leq
\infty$, and let $t_{0}$ be as in (\ref{T0}).
Proceeding similarly to the proof of Item (i), we have that: for  $\left\vert t\right\vert \in\lbrack t_{0}-\delta,t_{0}+\delta]$,
\begin{equation}
\sup_{\left\vert t\right\vert \in\lbrack t_{0}-\delta,t_{0}+\delta]}\left(
\max\{e^{(\alpha+h)|t|},|t|^{\widetilde{\alpha}+h}\}\left\Vert \mathcal{T}%
(u)-\mathcal{T}(v)\right\Vert _{(b+1,d)}\right)  \leq C_{4}\left\Vert
u-v\right\Vert _{E_{\alpha+h,\widetilde{\alpha}+h}^{d}}\left(  \left\Vert
u\right\Vert _{E_{\alpha,\widetilde{\alpha}}}^{b-1}+\left\Vert v\right\Vert
_{E_{\alpha,\widetilde{\alpha}}}^{b-1}\right)  , \label{ine888}%
\end{equation}
and for $t\geq t_{0}+\delta,$
\begin{align*}
&  \left\Vert \mathcal{T}(u)-\mathcal{T}(v)\right\Vert _{(b+1,d)}\\
&  \leq C\sup_{t>\delta}e^{(\alpha+h)t}\left\Vert (u-v)(t)\right\Vert
_{(b+1,d)}\left(  \sup_{t>\delta}e^{(b-1)\alpha t}\left\Vert u(t)\right\Vert
_{(b+1,\infty)}^{b-1}+\sup_{t>\delta}e^{(b-1)\alpha t}\left\Vert
v(t)\right\Vert _{(b+1,\infty)}^{b-1}\right) \\
&  \times\int_{t_{0}}^{t-\delta}s^{\frac{2}{b+1}}e^{-\beta s}e^{(-b\alpha
-h)(t-s)}ds\\
&  +C\sup_{0<t<\delta}t^{\widetilde{\alpha}+h}\left\Vert (u-v)(t)\right\Vert
_{(b+1,d)}\left(  \sup_{0<t<\delta}t^{(b-1)\widetilde{\alpha}}\left\Vert
u(t)\right\Vert _{(b+1,\infty)}^{b-1}+\sup_{0<t<\delta}t^{(b-1)\widetilde
{\alpha}}\left\Vert v(t)\right\Vert _{(b+1,\infty)}^{b-1}\right) \\
&  \times\int_{t-\delta}^{t}s^{\frac{2}{b+1}}e^{-\beta s}(t-s)^{(-b\widetilde
{\alpha}-h)}ds\\
&  +C\sup_{t>\delta}e^{(\alpha+h)t}\left\Vert (u-v)(t)\right\Vert
_{(b+1,d)}\left(  \sup_{t>\delta}e^{(b-1)\alpha t}\left\Vert u(t)\right\Vert
_{(b+1,\infty)}^{b-1}+\sup_{t>\delta}e^{(b-1)\alpha t}\left\Vert
v(t)\right\Vert _{(b+1,\infty)}^{b-1}\right) \\
&  \times\int_{0}^{t_{0}}s^{-\beta}e^{(-b\alpha-h)(t-s)}ds.
\end{align*}
Therefore, for $|t|\geq t_{0}+\delta$, we obtain that
\begin{equation}
\sup_{|t|\geq t_{0}+\delta}e^{(\alpha+h)|t|}\left\Vert \mathcal{T}%
(u)-\mathcal{T}(v)\right\Vert _{(b+1,d)}\leq C_{5}\left\Vert u-v\right\Vert
_{E_{\alpha+h,\widetilde{\alpha}+h}^{d}}\left(  \left\Vert u\right\Vert
_{E_{\alpha,\widetilde{\alpha}}}^{b-1}+\left\Vert v\right\Vert _{E_{\alpha
,\widetilde{\alpha}}}^{b-1}\right)  . \label{ine8}%
\end{equation}
Similarly, for $|t|<t_{0}$, it follows that
\begin{equation}
\sup_{|t|<t_{0}}|t|^{\widetilde{\alpha}+h}\left\Vert \mathcal{T}%
(u)-\mathcal{T}(v)\right\Vert _{(b+1,d)}\leq C_{6}\sup_{|t|<t_{0}%
}|t|^{\widetilde{\alpha}+h}\left\Vert u-v\right\Vert _{(b+1,d)}\left(
|t|^{(b-1)\widetilde{\alpha}}\left\Vert u\right\Vert _{(b+1,\infty)}%
^{b-1}+|t|^{(b-1)\widetilde{\alpha}}\left\Vert v\right\Vert _{(b+1,\infty
)}^{b-1}\right) . \label{ine88}%
\end{equation}

Putting together estimates (\ref{ine888}), (\ref{ine8}) and (\ref{ine88})
considering $\mathcal{K}_{d,h}=C_{4}+C_{5}+C_{6},$ and recalling the space
$E_{\alpha+h,\widetilde{\alpha}+h}^{d}$ (see (\ref{space-1})), we arrive at
\begin{equation}
\left\Vert \mathcal{T}(u)-\mathcal{T}(v)\right\Vert _{E_{\alpha+h,\widetilde
{\alpha}+h}^{d}}\leq\mathcal{K}_{d,h}\left\Vert u-v\right\Vert _{E_{\alpha
+h,\widetilde{\alpha}+h}^{d}}\left(  \left\Vert u\right\Vert _{E_{\alpha
,\widetilde{\alpha}}}^{b-1}+\left\Vert v\right\Vert _{E_{\alpha,\widetilde
{\alpha}}}^{b-1}\right)  . \label{ine8888}%
\end{equation}

Next the solution $u$ obtained in item (i) can be approximated by the Picard
sequence $\left\{  v_{m}\right\}  _{m\in\mathbb{N}}$ defined as%
\begin{equation}
v_{1}=\dot{W}(t)u_{0}+W(t)u_{1}\text{ and }v_{m+1}=\dot{W}(t)u_{0}%
+W(t)u_{1}+\int_{0}^{t}W(t-s)v_{m}(s)ds,\label{Picard}%
\end{equation}
where the limit is taken in the norm $\left\Vert \cdot\right\Vert
_{E_{\alpha,\widetilde{\alpha}}}$. Also, by the proof of Item (i), we have
that%
\begin{equation}
\left\Vert v_{m}\right\Vert _{E_{\alpha,\widetilde{\alpha}}}\leq
2\varepsilon.\label{size-sol-1}%
\end{equation}
In the sequel we show that the sequence (\ref{Picard}) is uniformly bounded in
the norm $\left\Vert \cdot\right\Vert _{E_{\alpha+h,\widetilde{\alpha}+h}^{d}%
}$. For that, denote
\[
\Gamma_{m,h}^{d}=\left\Vert v_{m}\right\Vert _{E_{\alpha+h,\widetilde{\alpha
}+h}^{d}}.
\]
First, for the case $h=0$, applying inequality (\ref{ine8888}) with $h=0$ and
$v=0$ yields%
\begin{align}
\Gamma_{m+1,0}^{d} &  \leq\sup_{\left\vert t\right\vert \geq t_{0}}%
e^{\alpha\left\vert t\right\vert }\left\Vert \dot{W}(t)u_{0}+W(t)u_{1}%
\right\Vert _{(b+1,d)}+\sup_{t\in(-t_{0},t_{0})}\left\vert t\right\vert
^{\widetilde{\alpha}}\left\Vert \dot{W}(t)u_{0}+W(t)u_{1}\right\Vert
_{(b+1,d)}\nonumber\\
&+\mathcal{K}_{d,0}\left\Vert v_{m}\right\Vert_{E_{\alpha
,\widetilde{\alpha}}^{d}}\left\Vert v_{m}\right\Vert_{E_{\alpha
,\widetilde{\alpha}}}^{b-1}\nonumber\\
&  \leq\Gamma_{1,0}^{d}+\left(  K_{\alpha,\widetilde{\alpha},0}\left\Vert
v_{m}\right\Vert _{E_{\alpha,\widetilde{\alpha}}}^{b-1}\right)  \Gamma
_{m,0}^{d} \nonumber\\
&  \leq\frac{1}{1-L_{d}}\Gamma_{1,0}^{d},\label{aux-ine}%
\end{align}
provided that $\varepsilon>0$ satisfies%
\begin{equation}
L_{d}=\mathcal{K}_{d,0}2^{b-1}\varepsilon^{b-1}<1.\label{ine11}%
\end{equation}
For the case $h\neq0$, applying again (\ref{ine8888}), we can estimate%
\begin{align}
\Gamma_{m+1,h}^{d} &  \leq\sup_{\left\vert t\right\vert \geq t_{0}}%
e^{(\alpha+h)\left\vert t\right\vert }\left\Vert \dot{W}(t)u_{0}%
+W(t)u_{1}\right\Vert _{(b+1,d)}+\sup_{t\in(-t_{0},t_{0})}\left\vert
t\right\vert ^{(\widetilde{\alpha}+h)}\left\Vert \dot{W}(t)u_{0}%
+W(t)u_{1}\right\Vert _{(b+1,d)} \nonumber\\
&  +\mathcal{K}_{d,h}\left\Vert v_{m}\right\Vert _{E_{\alpha+h,\widetilde
{\alpha}+h}^{d}}\left\Vert v_{m}\right\Vert _{E_{\alpha,\widetilde{\alpha}}%
}^{b-1} \nonumber\\
&  \leq\Gamma_{1,h}^{d}+\left(  \mathcal{K}_{d,h}\left\Vert v_{m}\right\Vert
_{E_{\alpha,\widetilde{\alpha}}}^{b-1}\right)  \Gamma_{m,h}^{d},\label{ine12}%
\end{align}
for all $m\in\mathbb{N}$. Taking $\varepsilon>0$ such that $L_{d,h}%
=\mathcal{K}_{d,h}2^{b-1}\varepsilon^{b-1}<1$ and using (\ref{cond-d-h}), the
R.H.S. of (\ref{ine12}) can be bounded by $\frac{1}{1-L_{d,h}}\Gamma_{1,h}%
^{d}<\infty$ and then we obtain the desired boundedness. Finally, the uniform
boundedness of $\{v_{m}\}_{m\in\mathbb{N}}$ in $E_{\alpha+h,\widetilde{\alpha
}+h}^{d}$ and the uniqueness of the limit in the sense of distributions yield
the property (\ref{Globalregular}).\fin

\subsection{Exponential stability}

\begin{theorem}
\label{Stable}(Exponential stability) Under the same assumptions of Theorem
\ref{Global}. Suppose also that $u$ and $\tilde{u}$ belonging to
$E_{\alpha,\widetilde{\alpha}}$ are solutions of (\ref{intergralEq}) obtained
in Theorem \ref{Global} with initial data $(u_{0},u_{1})$ and $(\tilde{u}%
_{0},\tilde{u}_{1})$, respectively. Then, we have that
\begin{equation}
\lim_{|t|\rightarrow\infty}e^{(\alpha+h)|t|}\left\Vert \dot{W}(t)(u_{0}%
-\tilde{u}_{0})+W(t)(u_{1}-\tilde{u}_{1})\right\Vert _{(b+1,d)}%
=0\label{condition}%
\end{equation}
if and only if
\begin{equation}
\lim_{|t|\rightarrow\infty}e^{(\alpha+h)|t|}\left\Vert u(t)-\tilde
{u}(t)\right\Vert _{(b+1,d)}=0.\label{Stability}%
\end{equation}
In particular, condition (\ref{condition}) holds provided that $D(u_{0}%
-\tilde{u}_{0})\in L^{(\frac{b+1}{b},d)}(\mathbb{H}^{n})$ and $(u_{1}%
-\tilde{u}_{1})\in L^{(\frac{b+1}{b},d)}(\mathbb{H}^{n})$.
\end{theorem}

\textbf{Proof.} First assuming (\ref{condition}), we prove (\ref{Stability}).
Without loss of generality, assume also that $t>0$. Then, we have that%
\begin{align}
e^{(\alpha+h)t}\left\Vert u(t)-\tilde{u}(t)\right\Vert _{(b+1,d)}  &  \leq
e^{(\alpha+h)t}\left\Vert \dot{W}(t)(u_{0}-\tilde{u}_{0})+W(t)(u_{1}-\tilde
{u}_{1})\right\Vert _{(b+1,d)} \nonumber\\
&  +e^{(\alpha+h)t}\int_{0}^{t}\left\Vert W(t-s)(F(u)-F(\tilde{u}%
))(s)ds\right\Vert _{(b+1,d)}ds. \label{ine-13}%
\end{align}
Take $\varepsilon>0$ such that $L_{d,h}=\mathcal{K}_{d,h}2^{b-1}%
\varepsilon^{b-1}<1$, where $\mathcal{K}_{d,h}$ is given in the proof of
Theorem \ref{Global} (ii) (see (\ref{ine12})). Therefore, it follows that
\begin{align}
\sup_{t\geq t_{0}}e^{(\alpha+h)t}\left\Vert u(t)\right\Vert _{(b+1,d)}%
+\sup_{t\in(0,t_{0})}t^{\widetilde{\alpha}+h}\left\Vert u(t)\right\Vert
_{(b+1,d)}  &  \leq\frac{1}{1-L_{d,h}}\Gamma_{1,h}^{d},\label{aux-bound-1}\\
\sup_{t\geq t_{0}}e^{(\alpha+h)t}\left\Vert \tilde{u}(t)\right\Vert
_{(b+1,d)}+\sup_{t\in(0,t_{0})}t^{\widetilde{\alpha}+h}\left\Vert \tilde
{u}(t)\right\Vert _{(b+1,d)}  &  \leq\frac{1}{1-L_{d,h}}\Gamma_{1,h}^{d}.
\label{aux-bound-2}%
\end{align}
provided that $\left\Vert (u_{0},u_{1})\right\Vert _{\mathcal{I}_{0}%
},\,\left\Vert (\tilde{u}_{0},\tilde{u}_{1})\right\Vert _{\mathcal{I}_{0}%
}\leq \varepsilon$. For $t\geq t_{0}+\delta$, using $\left\Vert
u\right\Vert _{E_{\alpha,\widetilde{\alpha}}},\left\Vert \tilde{u}\right\Vert
_{E_{\alpha,\widetilde{\alpha}}}\leq2\varepsilon$, we can handle the R.H.S. of
(\ref{ine-13}) as follows
\begin{align}
&  e^{(\alpha+h)t}\int_{0}^{t}\left\Vert W(t-s)(F(u)-F(\tilde{u}%
))(s)ds\right\Vert _{(b+1,d)}ds \nonumber\\
&  \leq C2^{b}\varepsilon^{b-1}e^{(\alpha+h)t}\int_{0}^{t_{0}}s^{-\beta
}e^{(-b\alpha-h)(t-s)}e^{(\alpha+h)(t-s)}\left\Vert u(t-s)-\tilde
{u}(t-s)\right\Vert _{(b+1,d)}ds \nonumber\\
&  +C2^{b}\varepsilon^{b-1}e^{(\alpha+h)t}\int_{t_{0}}^{t-\delta} s^{\frac{2}{b+1}} e^{-\beta
s}e^{(-b\alpha-h)(t-s)}e^{(\alpha+h)(t-s)}\left\Vert u(t-s)-\tilde
{u}(t-s)\right\Vert _{(b+1,d)}ds \nonumber\\
&  +C2^{b}\varepsilon^{b-1}e^{(\alpha+h)t}\int_{t-\delta}^{t} s^{\frac{2}{b+1}} e^{-\beta
s}(t-s)^{-b\widetilde{\alpha}-h}(t-s)^{\widetilde{\alpha}+h}\left\Vert
u(t-s)-\tilde{u}(t-s)\right\Vert _{(b+1,d)}ds \nonumber\\
&  \leq C2^{b}\varepsilon^{b-1}e^{-(b-1)\alpha t}e^{(b\alpha+h)t_{0}}\int
_{0}^{t_{0}}s^{-\beta}ds\left(  \left\Vert u\right\Vert _{E_{\alpha
+h,\widetilde{\alpha}+h}}+\left\Vert \tilde{u}\right\Vert _{E_{\alpha
+h,\widetilde{\alpha}+h}}\right) \nonumber\\
&  +C2^{b}\varepsilon^{b-1}e^{-((\beta-\sigma)-\alpha-h)t_{0}}\int_{t_{0}}^{t-\delta
}e^{-(b-1)\alpha(t-s)}ds\left(  \left\Vert u\right\Vert _{E_{\alpha
+h,\widetilde{\alpha}+h}}+\left\Vert \tilde{u}\right\Vert _{E_{\alpha
+h,\widetilde{\alpha}+h}}\right) \hbox{   (because   } s^{\frac{2}{b+1}}e^{-\sigma s}\leq C) \nonumber\\
&  +C2^{b}\varepsilon^{b-1}e^{-((\beta-\sigma)-\alpha-h)t}e^{\beta\delta}\int
_{t-\delta}^{t}(t-s)^{-b\widetilde{\alpha}-h}ds\left(  \left\Vert u\right\Vert
_{E_{\alpha+h,\widetilde{\alpha}+h}}+\left\Vert \tilde{u}\right\Vert
_{E_{\alpha+h,\widetilde{\alpha}+h}}\right)  . \label{aux-asymp-1}%
\end{align}
where $t_{0}$ is as in (\ref{T0}), $\delta=t_{0}/2$, and $\beta=\frac{n-1}%
{2}\left(  1-\frac{2}{b+1}\right)  $. Combining inequalities (\ref{ine-13}),
(\ref{aux-asymp-1}), (\ref{aux-bound-1})-(\ref{aux-bound-2}), and condition
(\ref{condition}), and using $b>1$ and $\beta-\sigma>\alpha+h$, we obtain that
\begin{align*}
\lim\sup_{t\rightarrow\infty}e^{(\alpha+h)t}\left\Vert u(t)-\tilde
{u}(t)\right\Vert _{(b+1,d)}\leq &  \lim\sup_{t\rightarrow\infty}%
e^{(\alpha+h)t}\left\Vert \dot{W}(t)(u_{0}-\tilde{u}_{0})+W(t)(u_{1}-\tilde
{u}_{1})\right\Vert _{(b+1,d)}\\
&  +\lim\sup_{t\rightarrow\infty}e^{(\alpha+h)t}\int_{0}^{t}\left\Vert
W(t-s)(F(u)-F(\tilde{u}))(s)\right\Vert _{(b+1,d)}ds\\
&  =\lim\sup_{t\rightarrow\infty}e^{(\alpha+h)t}\left\Vert \dot{W}%
(t)(u_{0}-\tilde{u}_{0})+W(t)(u_{1}-\tilde{u}_{1})\right\Vert _{(b+1,d)}\\
&  +C2^{b}\varepsilon^{b-1}e^{(b\alpha+h)t_{0}}t_{0}^{1-\beta}\frac
{2}{1-L_{d,h}}\Gamma_{1,h}^{d}\text{ }e^{-(b-1)\alpha t}\\
&  +C2^{b}\varepsilon^{b-1}\frac{2}{1-L_{d,h}}\Gamma_{1,h}^{d}\text{
}e^{-(b-1)\alpha t}\left(  e^{-((\beta-\sigma)-b\alpha-h)t_{0}}-e^{-((\beta-\sigma)
-b\alpha-h)(t-\delta)}\right) \\
&  +C2^{b}\varepsilon^{b-1}\frac{2}{1-L_{d,h}}\Gamma_{1,h}^{d}\text{ }%
e^{\beta\delta}\text{ }\delta^{1-b\widetilde{\alpha}-h}e^{-((\beta-\sigma)-\alpha
-h)t}\\
&  \rightarrow0+0+0+0=0,
\end{align*}
which implies (\ref{Stability}). Next, for the reciprocal assertion, assume
(\ref{Stability}). Then, using the same bounds for
\[
e^{(\alpha+h)t}\int_{0}^{t}\left\Vert W(t-s)(F(u)-F(\tilde{u}%
))(s)ds\right\Vert _{(b+1,d)}ds,
\]
we arrive at%
\begin{align*}
&\lim\sup_{t\rightarrow\infty}e^{(\alpha+h)t}\left\Vert \dot{W}(t)(u_{0}%
-\tilde{u}_{0})+W(t)(u_{1}-\tilde{u}_{1})\right\Vert _{(b+1,d)} \\
&  \leq\lim\sup_{t\rightarrow\infty}e^{(\alpha+h)t}\left\Vert u(t)-\tilde
{u}(t)\right\Vert _{(b+1,d)}+\lim\sup_{t\rightarrow\infty}e^{(\alpha+h)t}\int_{0}^{t}\left\Vert
W(t-s)(F(u)-F(\tilde{u}))(s)\right\Vert _{(b+1,d)}ds\\
&  =0+0=0,
\end{align*}
which gives the desired conclusion. \fin

\subsection{Scattering and wave operator}

In this part we analyze the scattering property of the global mild
solution obtained in Theorem \ref{Global}.

\begin{theorem}
\label{SScattering}(Scattering) Suppose the same conditions of Theorem
\ref{Global}. Consider also the solution $u$ of (\ref{intergralEq}) obtained in Theorem
\ref{Global} with initial data $(u_{0},u_{1})$. Then, there exist $(u_{0}%
^{\pm},u_{1}^{\pm})\in\mathcal{I}_{0}$ satisfying
\begin{equation}
\left\Vert u(t)-u^{+}(t)\right\Vert _{(b+1,d)}=O(e^{-b(\alpha+h)t}),\hbox{
as  }t\rightarrow+\infty,\label{ine16}%
\end{equation}%
\begin{equation}
\left\Vert u(t)-u^{-}(t)\right\Vert _{(b+1,d)}=O(e^{-b(\alpha+h)|t|}%
),\hbox{  as  }t\rightarrow-\infty,\label{ine17}%
\end{equation}
where $u^{\pm}$ are the unique solutions of the associated linear problem
$u^{\pm}=\dot{W}(t)u_{0}^{\pm}+W(t)u_{1}^{\pm}.$
\end{theorem}

\textbf{Proof.} We show only the property (\ref{ine16}). The proof of
(\ref{ine17}) is left to the reader. We start by defining
\[
u_{1}^{+}=u_{1}+\int_{0}^{\infty}W(-s)F(u(s))ds\text{ and }\,u_{0}^{+}=u_{0}.
\]
For $t>0$, consider $u^{+}$ given by
\[
u^{+}=\dot{W}(t)u_{0}+W(t)u_{1}+\int_{0}^{\infty}W(t-s)F(u(s))ds,
\]
and note that
\begin{equation}
u-u^{+}=\int_{t}^{\infty}W(s-t)F(u(s))ds=\int_{0}^{\infty}W(s)F(u(t+s))ds.
\label{Scat1}%
\end{equation}
Using (\ref{Scat1}) with $t>t_{0}$ and recalling $0<\beta=\frac{n-1}{2}\left(
1-\frac{2}{b+1}\right)  <1,$ we have that
\begin{align*}
&  \left\Vert u-u^{+}\right\Vert _{(b+1,d)}=\left\Vert \int_{0}^{\infty
}W(s)F(u(t+s))ds\right\Vert _{(b+1,d)}\\
&  \leq\left(  \int_{0}^{t_{0}}s^{-\beta}e^{-b(\alpha+h)(s+t)}ds+\int_{t_{0}}^{\infty}s^{\frac{2}{b+1}}e^{-\beta s}e^{-b(\alpha+h)(s+t)}ds\right) \\
&  \times\left(  \sup_{t\geq t_{0}}e^{(\alpha+h)t}\left\Vert u(t)\right\Vert
_{(b+1,d)}\right)^b\\
&  \leq C\left(  \Gamma_{1,h}^{d}\right)  ^{b}\left(
\int_{0}^{t_{0}}s^{-\beta}e^{-b(\alpha+h)(s+t)}ds+\int_{t_{0}}^{\infty}s^{\frac{2}{b+1}}
e^{(-1-\sigma-\alpha-bh)(s+t)}ds\right) \\
&\hbox{(because   } 1-\beta + \sigma - (b-1)\alpha=0\Rightarrow \beta+b(\alpha +h) = 1+\sigma+\alpha+bh)\\
&  \leq C\left(  \Gamma_{1,h}^{d}\right)  ^{b}\left(e^{-b(\alpha+h)t}
\int_{0}^{t_{0}}s^{-\beta}e^{-b(\alpha+h)s}ds + e^{(-\beta-b(\alpha+h))t}\int_{t_{0}}^{\infty}e^{(-1-\alpha-bh)s}ds\right) \\
&\hbox{(because   } s^{\frac{2}{b+1}}e^{-\sigma s} <C)\\
&  \leq C\left(  \Gamma_{1,h}^{d}\right)  ^{b}e^{-b(\alpha+h)t} \left(
\int_{0}^{t_{0}}s^{-\beta} ds + \int_{t_{0}}^{\infty}e^{(-1-\alpha-bh)s}ds\right) \\
&  \leq C\left(\frac{t_0^{1-\beta}}{1-\beta}+ \frac{e^{(1+\alpha+bh)t_0}}{1+\alpha+bh} \right)\left(  \Gamma_{1,h}^{d}\right)  ^{b}e^{-b(\alpha+h)t}.
\end{align*}
This leads to the scattering behaviour \eqref{ine16}. Our proof is completed. \fin

\bigskip

\begin{remark}
\label{scat-2} In view of Theorem \ref{Stable}, we can improve the scattering
decay by replacing $O$ with $o$ in Theorem \ref{SScattering}. Precisely,
letting $u^{\pm}$ be as in Theorem \ref{SScattering}, we have that
\begin{equation}
\left\Vert u(t)-u^{\pm}(t)\right\Vert _{(b+1,d)}=o(e^{-b(\alpha+h)|t|}),\text{
as }t\rightarrow\pm\infty, \label{improve2}%
\end{equation}
provided that
\begin{equation}
\lim_{t\rightarrow\pm\infty}e^{(\alpha+h)|t|}\left\Vert \dot{W}(t)u_{0}%
+W(t)u_{1}\right\Vert _{(b+1,d)}=0. \label{improve1}%
\end{equation}
In fact, by (\ref{improve1}), Theorem \ref{Stable} gives
\begin{equation}
\lim_{t\rightarrow\pm\infty}e^{(\alpha+h)\left\vert t\right\vert }\left\Vert
u(t)\right\Vert _{(b+1,d)}=0. \label{improve3}%
\end{equation}
It follows that
\begin{align*}
&  \left\Vert u(t)-u^{+}(t)\right\Vert _{(b+1,d)}\\
&  =\left\Vert \int_{0}^{\infty}W(s)F(u(s+t))ds\right\Vert _{(b+1,d)}\\
&  \leq2^{b-1}(\Gamma_{1,h}^{d})^{b-1}\left(  e^{-b(\alpha+h)t}\int_{0}%
^{t_{0}}s^{-\beta}e^{-b(\alpha+h)s}e^{(\alpha+h)(s+t)}\left\Vert
u(s+t)\right\Vert _{(b+1,d)}ds\right. \\
&\left.+ e^{(-\beta-b(\alpha+h))t}\int_{t_{0}}^{\infty}s^{\frac{2}{b+1}}e^{-(1-\sigma-\alpha-bh)s}e^{(\alpha+h)(s+t)}\left\Vert u(s+t)\right\Vert _{(b+1,d)}ds\right)\\
&  \leq C2^{b-1}(\Gamma_{1,h}^{d})^{b-1}\left(  e^{-b(\alpha+h)t}\int_{0}%
^{t_{0}}s^{-\beta}e^{-b(\alpha+h)s}e^{(\alpha+h)(s+t)}\left\Vert
u(s+t)\right\Vert _{(b+1,d)}ds\right. \\
&\left.+ e^{(-\beta-b(\alpha+h))t}\int_{t_{0}}^{\infty} e^{-(1-\alpha-bh)s}e^{(\alpha+h)(s+t)}\left\Vert u(s+t)\right\Vert _{(b+1,d)}ds\right) \hbox{   (because   } s^{\frac{2}{b+1}}e^{-\sigma s}<C),
\end{align*}
and then
\begin{align*}
\limsup_{t\rightarrow\infty}e^{b(\alpha+h)t}\left\Vert u(t)-u^{+}%
(t)\right\Vert _{(b+1,h)} &\leq C2^{b-1}(\Gamma_{1,h}^{d})^{b-1}\left(  \int_{0}^{t_{0}}s^{-\beta
}ds + \int_{t_{0}}^{\infty}e^{-(1+\alpha+bh)s}ds\right) \\
&  \times\limsup_{t\rightarrow\infty}e^{(\alpha+h)t}\left\Vert u(t)\right\Vert
_{(b+1,d)}\\
&\leq C2^{b-1}(\Gamma_{1,h}^{d})^{b-1}\left(  \frac{t_0^{1-\beta}}{1-\beta} + \frac{e^{-(1+\alpha+bh)t_0}}{1+\alpha+bh}\right) \\
&  \times\limsup_{t\rightarrow\infty}e^{(\alpha+h)t}\left\Vert u(t)\right\Vert
_{(b+1,d)}\rightarrow0,
\end{align*}
as desired. The case $t\rightarrow-\infty$ is done by a similar way.
\end{remark}

Finally, we establish the construction of wave operators for equation \eqref{Klein1}, that is, the construction of solutions with prescribed scattering states.  In particular, for a given initial data $(f_0,f_1)$, we find a global solution $u$ which converges as $|t| \to \infty$ to the solution of the associated linear problem with initial data $(f_0,f_1)$. In the following theorem we state and prove the construction of future wave operator, the past wave one is constructed by the same way.

\begin{theorem}\label{WaveOp}(Future wave operator)
Suppose the same conditions of Theorem \ref{Global} and let $\gamma \in [0,\infty)$. Then, for any $(f_0,f_1)\in \mathcal{S}'(\mathbb{H}^n)\times \mathcal{S}'(\mathbb{H}^n)$ such that $(Df_0,f_1) \in L^{\left(\frac{b+1}{b},\infty \right)}(\mathbb{H}^n)\times L^{\left(\frac{b+1}{b},\infty \right)}(\mathbb{H}^n)$ and
$$\sup_{t\geq t_0}e^{\gamma t}\norm{\dot{W}(t)f_0 + W(t)f_1}_{(b+1,\infty)} <\infty,$$
there exist $T_0 = T_0(f_0,f_1)>t_0$ and a solution $u$ of the integral equation \eqref{intergralEq} on $[T_0,\infty)$ satisfying
$$\sup_{t>T_0} e^{{\gamma}t}\norm{u(t)}_{(b+1,\infty)} <\infty$$
and
$$\lim_{t\to \infty} e^{\gamma t}\norm{u(t) - \left(\dot{W}(t)f_0 + W(t)f_1 \right)}_{(b+1,\infty)}=0.$$
If $v(t)$ is also a solution of equation \eqref{intergralEq} on $[T_0,\infty)$ satisfying
$$\sup_{t>T_0} e^{\gamma t}\norm{v(t)}_{(b+1,\infty)}<\infty,$$
then there exists $T>T_0$ such that $u=v$ on $[T,\infty)$.
\end{theorem}
\begin{proof}
For $T>t_0$, we consider the space $E_{\gamma}^{\geq T}$ of all Bochner measurable
$u:[T,\infty)\rightarrow L^{(b+1,\infty)}(\mathbb{H}^n)$ satisfying
\begin{equation}
\left\Vert u\right\Vert _{E_{\gamma}^{\geq T}}=\sup_{t \geq T}e^{\gamma t}\left\Vert
u(t)\right\Vert _{(b+1,\infty)} < \infty
.\label{space-1}%
\end{equation}
The space $E_{\gamma}^{\geq T}$ endowed with $\left\Vert
\cdot\right\Vert _{E_{\gamma}^{\geq T}}$ is a Banach space.

\underline{Step 1: Construction of the solution.} Let $R>0$ and consider the closed ball
$$B^{\geq T}(0,R) = \left\{ u \in  E_{\gamma}^{\geq T}: \norm{u}_{E_{\gamma}^{\geq T}} \leq R  \right\}.$$
We define the mapping $\Phi: B^{\geq T}(0,R) \to B^{\geq T}(0,R)$ by
\begin{equation*}
\Phi(\omega) = -\int_t^\infty W(t-s) F\left(\omega(s) + \dot{W}(s)f_0 + W(s)f_1 \right) ds.
\end{equation*}
Observe that, if $\omega$ is a fixed point of the operator $\Phi: B^{\geq T_0}(0,R) \to B^{\geq T_0}(0,R)$, for some $T_0>0$ which will be chosen later, then
$$u(t)= \dot{W}(t)f_0 + W(t)f_1 + \omega(t)$$
satisfies $u \in E_\gamma^{\geq T_0}$ and $u\in C([T_0,\infty); \mathcal{S}'(\mathbb{H}^n))$. Using the group properties of $W(t)$, we can show that
\begin{equation}\label{wave1}
u(t) = \dot{W}(t-T_0)u(T_0) + W(t-T_0)\partial_t u(T_0) + \int_{T_0}^tW(t-s)F(u(s))ds,
\end{equation}
where $\mathcal{G}^{t\geq T_0}(u) = \int_{T_0}^tW(t-s)F(u(s))ds \in \mathcal{S}'(\mathbb{H}^n)$ and equation \eqref{wave1} holds for a.e. $t\geq T_0$. The function $u$ given by \eqref{wave1} is a solution of the integral equation \eqref{Klein1} on $[T_0,\infty)$.

\underline{Step 2: The existence of unique fixed point $\omega$ of $\Phi$.} We prove this argument by establishing that $\Phi$ is a contraction mapping. Indeed, by using (\ref{T0}), Lemma \ref{DisperInter}, Remark \ref{Re} (ii) and H\"{o}lder's inequality, for $t\geq T>t_0$, we estimate
\begin{align*}
&\norm{\Phi(\omega)}_{(b+1,\infty)} = \left\Vert \int_t^{\infty
}W(s-t)F \left( \omega(s) + \dot{W}(s)f_0 + W(s)f_1 \right)ds\right\Vert _{(b+1,\infty)}\nonumber\\
&=\left\Vert \int_{0}^{\infty
}W(s)F \left( \omega(t+s) + \dot{W}(t+s)f_0 + W(t+s)f_1 \right)ds\right\Vert _{(b+1,\infty)}\nonumber\\
&  \leq \left(  \int_{0}^{t_{0}}s^{-\frac{n-1}{2}\left(  1-\frac{2}%
{b+1}\right)  }e^{-b\gamma(s+t)}ds+\int_{t_{0}}^{\infty}s^{\frac{2}{b+1}}e^{-\frac{n-1}%
{2}\left(  1-\frac{2}{b+1}\right)  s}e^{-b\gamma(s+t)}ds\right) \\
&\times\left(  \sup_{t\geq T}e^{\gamma t}\left\Vert \omega(t) + \dot{W}(t)f_0 + W(t)f_1 \right\Vert
_{(b+1,\infty)}\right)^{b}\\
&  \leq \left(  \sup_{t\geq T}e^{\gamma t}\norm{\omega(t)}_{(b+1,\infty)} + \sup_{t\geq T} e^{\gamma t}\norm{\dot{W}(t)f_0 + W(t)f_1}_{(b+1,\infty)} \right)^{b}\\
&\times \left(
e^{-b\gamma t}\int_{0}^{t_{0}}s^{-\beta}e^{-b\gamma s}ds+ e^{-b\gamma t}\int_{t_{0}}^{\infty}%
(s^{\frac{2}{b+1}}e^{-\sigma s})e^{-(\beta-\sigma+b\gamma)s}ds\right) \\
&\leq C\left(  \sup_{t\geq T}e^{\gamma t}\norm{\omega(t)}_{(b+1,\infty)} + \sup_{t\geq T} e^{\gamma t}\norm{\dot{W}(t)f_0 + W(t)f_1}_{(b+1,\infty)} \right)^{b}\\
&\times e^{-\gamma t}e^{-(b-1)\gamma T}\left(
\int_{0}^{t_{0}}s^{-\beta} ds + \int_{t_{0}}^{\infty}e^{-b\gamma s}ds\right) \hbox{   (because   } \beta>\sigma \hbox{  and  } s^{\frac{2}{b+1}}e^{-\sigma s}<C)\\
&  \leq \left(  \sup_{t\geq T}e^{\gamma t}\norm{\omega(t)}_{(b+1,\infty)} + \sup_{t\geq T} e^{\gamma t}\norm{\dot{W}(t)f_0 + W(t)f_1}_{(b+1,\infty)} \right)^{b}e^{-\gamma t}C(T),
\end{align*}
where
\begin{align*}
C(T) &= e^{-(b-1)\gamma T}\left(
\int_{0}^{t_{0}}s^{-\beta} ds + \int_{t_{0}}^{\infty}e^{-b\gamma s}ds\right) \\
&= e^{-(b-1)\gamma T}\left(
\frac{t_0^{1-\beta}}{1-\beta} + \frac{e^{b\gamma t_0}}{b\gamma}\right) \to 0 \hbox{    as   } T\to \infty.
\end{align*}
Hence
$$\sup_{t\geq T}e^{\gamma t}\norm{\Phi(\omega)}_{(b+1,\infty)} \leq C(T)\left(  \sup_{t\geq T}e^{\gamma t}\norm{\omega(t)}_{(b+1,\infty)} + \sup_{t\geq T} e^{\gamma t}\norm{\dot{W}(t)f_0 + W(t)f_1}_{(b+1,\infty)} \right)^{b},$$
for $t\geq T\geq t_0$.
This leads to
\begin{equation}\label{core1}
\norm{\Phi(\omega)}_{E_\gamma^{\geq T}} \leq C(T)\left(  R + \norm{\dot{W}(t)f_0 + W(t)f_1}_{E_\gamma^{\geq T}} \right)^{b},
\end{equation}
for all $t\geq T\geq t_0$ and $\omega\in B^{\geq T}(0,R)$. By the same way and using the property \eqref{nonlinearity-1} of function $F$, for every $\omega,\, \widetilde{\omega}\in B^{\geq T}(0,R)$, we can estimate
\begin{align}\label{core2}
\sup_{t\geq T}e^{\gamma t}\norm{\Phi(\omega)-\Phi(\widetilde\omega)}_{(b+1,\infty)} &\leq {C}(T)\norm{\omega-\widetilde\omega}_{E_\gamma^{\geq T}}  \left[\left(\norm{\omega}_{E_\gamma^{\geq T}} + \norm{\dot{W}(t)f_0 + W(t)f_1}_{E_\gamma^{\geq T}} \right)^{b-1} \right.\nonumber\\
&\left. + \left(\norm{\omega}_{E_\gamma^{\geq T}} + \norm{\dot{W}(t)f_0 + W(t)f_1}_{E_\gamma^{\geq T}} \right)^{b-1} \right]\nonumber\\
&\leq 2{C}(T)\norm{\omega-\widetilde\omega}_{E_\gamma^{\geq T}}  \left(R + \norm{\dot{W}(t)f_0 + W(t)f_1}_{E_\gamma^{\geq T}} \right)^{b-1}.
\end{align}
From inequalities \eqref{core1} and \eqref{core2}, we can choose $T_0=T_0(f_0,f_1)\geq t_0$ such that $\Phi$ is a contraction mapping on $B^{\geq T_0}(0,R)$, which implies the existence of a unique fixed point $\omega$ of $\Phi$.

\underline{Step 3: The convergence and uniqueness.}
Let $0<\eta<T_0$ be fixed and $\tau>2T_0$. Hence, we have $\widetilde{\tau} = \tau-\eta>T_0$. By the same estimates in Step 2, we obtain that \begin{align*}
e^{\gamma \tau}\norm{u(\tau)-\left(\dot{W}(\tau)f_1 + W(\tau)f_0 \right)}_{(b+1,\infty)} &= e^{\gamma \tau}\norm{\omega(\tau)}_{(b+1,\infty)}\\
&\leq \sup_{t>\widetilde{\tau}} e^{\gamma t}\norm{\Phi(\omega)(t)}_{(b+1,\infty)}\\
&\leq C(\widetilde{\tau})\left(R+ \sup_{t>T_0}e^{\gamma t}\norm{\dot{W}(t)f_0 + W(t)f_1}_{(b+1,\infty)}\right)^b\\
&\to 0, \hbox{   as   } \widetilde{\tau}\to \infty,
\end{align*}
due to $C(\widetilde{\tau})\to 0$ as $\widetilde{\tau}\to \infty$.
Finally, the uniqueness assertion follows from standard arguments (see, for example, \cite{FaFe}).

\end{proof}

\vspace{0.3cm}

\noindent\textbf{Conflict of interest statement:} The authors declare that
they have no conflict of interest.


\begin{thebibliography}{99}                                                                                               %


\bibitem {Fe2017}M.F. de Almeida, L.C.F. Ferreira, \textit{Time-weighted
estimates in Lorentz spaces and self-similarity for wave equations with
singular potentials}, Analysis \& PDE 10 (2) (2017), 423-438.

\bibitem {Anker2012}J.-P. Anker, V. Pierfelice, M. Vallarino, \textit{The wave
equation on hyperbolic spaces}, J. Differential Equations 252 (2012), 5613-5661.

\bibitem {Anker2014}J.-P. Anker, V. Pierfelice, \textit{Wave and Klein-Gordon
equations on hyperbolic spaces}, Analysis \& PDE, Vol. 7, No. 4 (2014), 953-995.

\bibitem {An2000}F. Antoci, \textit{Scattering theory for $p$-forms on
hyperbolic real space}, Istituto Lombardo (Rend. Sc.) A 134, 71-85 (2000).

\bibitem {Ba2008}V. Banica, R. Carles, G. Staffilani, \textit{Scattering
Theory for Radial Nonlinear Schr\"{o}dinger Equations on Hyperbolic Space},
Geometric and Functional Analysis \textbf{18}, 367-399 (2008)

\bibitem {Ba2009}V. Banica, R. Carles, T. Duyckaerts, \textit{On scattering
for NLS: from Euclidean to hyperbolic space}, Discrete Contin. Dyn. Syst. 24
(4) (2009), 1113-1127.

\bibitem {Bae}J. C. Baez, I. E. Segal, Z.-F. Zhou, \textit{The global Goursat
problem and scattering for nonlinear wave equations}, J. Funct. Anal. 93
(1990), 239-269

\bibitem {BeLo}J. Bergh, J. L\"{o}fstr\"{o}m, \textit{Interpolation Spaces. An
introduction}, Grundlehren der mathematischen Wissenschaften, Springer, Berlin (1976).

\bibitem {Bel}E. Belchev, M. Kepka, Z. Zhou, \textit{Global existence of
solutions to nonlinear wave equations}, Comm. Partial Differential Equations
24 (1999), 2297-2331.

\bibitem {Bo2008}J.-F. Bony, D. H\"{a}fner, \textit{The Semilinear Wave
Equation on Asymptotically Euclidean Manifolds}, Communications in Partial
Differential Equations 35(1), 23-67 (2009).

\bibitem {Ca1998}T. Cazenave, F.B. Weissler, \textit{Asymptotically
self-similar global solutions of the nonlinear Schr\"{o}dinger and heat
equations}, Math. Z. 228 (1) (1998), 83-120.

\bibitem {Ca2000}T. Cazenave, F.B. Weissler, \textit{Scattering theory and
self-similar solutions for the nonlinear Schr\"{o}dinger equation}, SIAM J.
Math. Anal. 31 (3) (2000), 625-650.

\bibitem {Ca2001}T. Cazenave, L. Vega, M.C. Vilela, \textit{A note on the
nonlinear Schr\"{o}dinger equation in weak-$L^{p}$ space}, Communications in
Contemporary Mathematics  3 (1) (2001), 153-162.

\bibitem {De1992}S. Debi\`{e}vre, H.D. Hislop, I.M. Sigal, \textit{Scattering
theory for the wave equation on non-compact manifolds}, Reviews in
Mathematical Physics 4 (4) (1992), 575-618.

\bibitem{FaFe} L.G. Farah and L.C.F Ferreira, \textit{ On the wave operator for the generalized Boussinesq equation}, Proc. Amer. Math. Soc. 140, 3055–3066 (2012)

\bibitem {Fe2009}L.C.F. Ferreira, E.J. Vilamizar-Roa, P.B. E Silva,
\textit{On the existence of infinite energy solutions for nonlinear
Schr\"{o}dinger equations}, Proc. Amer. Math. Soc. 137 (6) (2009), 1977-1987.

\bibitem {Fe2011}L.C.F. Ferreira, \textit{Existence and scattering theory for
Boussinesq type equations with singular data}, J. Differential Equation 250
(5) (2011), 2372-2388.

\bibitem {Fec1982}H. Pecher, \textit{Decay of solutions of nonlinear wave
equations in three space dimensions}, J. Funct. Anal. 46 (1982), 221-229.

\bibitem {Fec1984}H. Pecher, \textit{Nonlinear small data scattering for the
wave and Klein-Gordon equation}, Math. Z. 185 (1984), 261-270.

\bibitem {Fr2012}A. French, \textit{Scattering For Nonlinear Waves On
Hyperbolic Space}, PhD's thesis, University of North Carolina at Chapel Hill
(2012), https://doi.org/10.17615/6fp0-qw21

\bibitem {Fri2001}F.G. Friedlander, \textit{Notes on the Wave Equation on
Asymptotically Euclidean Manifolds}, J. Funct. Anal. 184 (1) (2001), 1-18.

\bibitem {Ge1997}V. Georgiev, H. Lindblad, C. D. Sogge, \textit{Weighted
Strichartz estimates and global existence for semilinear wave equations},
Amer. J. Math. 119 (1997), 1291-1319.

\bibitem {Gi1}J. Ginibre, G. Velo, \textit{The global Cauchy problem for the
non linear Klein-Gordon equation}, Math. Z. 189 (1985), 487-505.

\bibitem {Gi2}J. Ginibre, G. Velo, \textit{Conformal invariance and time decay
for non linear wave equations}, I, Ann. Inst. H. Poincar\'{e} Phys. Th\'{e}or.
47 (1987), 221-261.

\bibitem {G3}J. Ginibre, G. Velo, \textit{Conformal invariance and time decay
for nonlinear wave equations}, II, Ann. Inst. H. Poincar\'{e} Phys. Th\'{e}or.
47 (1987), 263-276.

\bibitem {G4}J. Ginibre, G. Velo, \textit{Scattering theory in the energy
space for a class of non-linear wave equations}, Comm. Math. Phys. 123 (1989), 535-573.

\bibitem {Gl}R. T. Glassey, \textit{Finite-time blow-up for solutions of
nonlinear wave equations}, Math. Z. 177 (1981), 323-340.

\bibitem {Hi1998}K. Hidano, \textit{Nonlinear small data scattering for the
wave equation in $\mathbb{R}^{4+1}$}, J. Math. Soc. Japan 50 (1998), 253-292.

\bibitem {Hi2000}K. Hidano, \textit{Small data scattering and blow-up for a
wave equation with a cubic convolution}, Funkcial. Ekvac. 43 (2000), 559-588

\bibitem {Hi2001}K. Hidano, \textit{Scattering problem for the nonlinear wave
equation in the finite energy and conformal charged}, J. Funct. Anal. 187 (2)
(2001), 274-307.

\bibitem {Tao}M. Keel, T. Tao, \textit{Endpoint Strichartz estimates}, Amer.
J. Math. 120 (1998), 955-980.

\bibitem {Jo}R. T. Glassey, \textit{Finite-time blow-up for solutions of
nonlinear wave equations}, Math. Z. 177 (1981), 323-340.

\bibitem {Lax1967}P.D. Lax, R.S. Phillips, \textit{Scattering theory}, Pure
and Applied Mathematics, Vol. 26, Academic Press, New York-London, 1967, xii +
276 pp. (1 plate).

\bibitem {Lin1995}H. Lindblad, C. D. Sogge, \textit{On existence and
scattering with minimal regularity for semilinear wave equations}, J. Funct.
Anal. 130 (1995), 357-426.

\bibitem{Wang2014} H. Lindblad, J. Metcalfe,C.D. Sogge, M. Tohaneanu and Ch. Wang, {\it The Strauss conjecture on Kerr black hole backgrounds},  Math. Ann. 359, 637–661 (2014).

\bibitem {LiYu}T. T. Li, X. Yu, \textit{Life-span of classical solutions to
fully nonlinear wave equations}, Comm. Partial Differential Equations 16
(1991), 909-940.

\bibitem {LiZho}T. T. Li, Y. Zhou, \textit{A note on the life-span of
classical solutions to nonlinear wave equations in four space dimensions},
Indiana Univ. Math. J. 44 (1995), 1207-1248

\bibitem {Liu2009}S. Liu, \textit{Remarks on infinite energy solutions of
nonlinear wave equations}, Nonlinear Analysis 71 (2009) 4231-4240.

\bibitem {MeTa2011}J. Metcalfe, M. Taylor, \textit{Nonlinear waves on 3D
hyperbolic space},  Trans. Amer. Math. Soc. 363 (7) (2011), 3489-3529.

\bibitem {MeTa2012}J. Metcalfe, M. Taylor, \textit{Dispersive wave estimates
on 3D hyperbolic space},  Proc. Amer. Math. Soc. 140 (11) (2012), 3861-3866.

\bibitem{Me2017} J. Metcalfe and Ch. Wang, {\it The Strauss Conjecture on Asymptotically Flat Space-Times}, SIAM Journal on Mathematical Analysis, Vol. 49, Iss. 6 (2017)

\bibitem {Phi1987}R. Phillips, B. Wiskott, A. Woo, \textit{Scattering theory
for the wave equation on a hyperbolic manifold}, J. Funct. Anal. 74 (2)
(1987), 346-398.

\bibitem {Sa2005}A. S\'{a} Barreto, \textit{Radiation fields, scattering and
inverse scattering on asymptotically hyperbolic manifolds}, Duke Math. J. 129
(3) (2005), 407-480.

\bibitem {Sa2008}A. S\'{a} Barreto, \textit{A support theorem for the
radiation fields on asymptotically Euclidean manifolds}, Math. Res. Lett. 15
(5) (2008), 973-991.

\bibitem {Si}T. C. Sideris, \textit{Nonexistence of global solutions to
semilinear wave equations in high dimensions}, J. Differential Equations 52
(1984), 378-406

\bibitem {Str1}W. A. Strauss, \textit{Decay and Asymptotics for $\Box u=F(u)$%
}, J. Funct. Anal. 2 (1968), 409-457.

\bibitem {Str2}W. A. Strauss, \textit{Nonlinear scattering theory at low
energy}, J. Funct. Anal. 41 (1981), 110-133.

\bibitem{Wang19} Y. Sire, C. D. Sogge, C. Wang, {\it The Strauss conjecture on negatively curved backgrounds}, Discrete and Continuous Dynamical Systems, 2019, 39(12): 7081-7099

\bibitem{Wang2019} Ch. Wang and X. Zhang {\it Wave equations with logarithmic nonlinearity on hyperbolic spaces}, preprint (2023) https://doi.org/10.48550/arXiv.2304.01595


\bibitem {Ta2001}D. Tataru, \textit{Strichartz estimates in the hyperbolic
space and global existence for the semilinear wave equation},  Trans. Amer.
Math. Soc. 353 (2) (2001), 795-807.

\bibitem {Zha2015}J. Zhang, \textit{Strichartz estimates and nonlinear wave
equation on nontrapping asymptotically conic manifolds}, Adv. Math. 271
(2015), 91-111.

\bibitem {Zho1995}Y. Zhou, \textit{Cauchy problem for semilinear wave
equations in four space dimensions with small initial data}, J. Partial
Differential Equations 8 (2) (1995), 135-144.
\end{thebibliography}
\end{document}